\documentclass[11pt]{amsart}
\usepackage{amsbsy}
\usepackage{graphicx,epsfig,subfigure,psfrag}
\begin{document}

\newtheorem{theorem}{Theorem}
\newtheorem{proposition}{Proposition}
\newtheorem{lemma}{Lemma}
\newtheorem{corollary}{Corollary}
\newtheorem{definition}{Definition}
\newtheorem{remark}{Remark}
\newcommand{\tex}{\textstyle}
\numberwithin{equation}{section} \numberwithin{theorem}{section}
\numberwithin{proposition}{section} \numberwithin{lemma}{section}
\numberwithin{corollary}{section}
\numberwithin{definition}{section} \numberwithin{remark}{section}
\newcommand{\ren}{\mathbb{R}^N}
\newcommand{\re}{\mathbb{R}}
\newcommand{\n}{\nabla}
\newcommand{\p}{\partial}
\newcommand{\iy}{\infty}
\newcommand{\pa}{\partial}
\newcommand{\fp}{\noindent}
\newcommand{\ms}{\medskip\vskip-.1cm}
\newcommand{\mpb}{\medskip}
\newcommand{\AAA}{{\bf A}}
\newcommand{\BB}{{\bf B}}
\newcommand{\CC}{{\bf C}}
\newcommand{\DD}{{\bf D}}
\newcommand{\EE}{{\bf E}}
\newcommand{\FF}{{\bf F}}
\newcommand{\GG}{{\bf G}}
\newcommand{\oo}{{\mathbf \omega}}
\newcommand{\Am}{{\bf A}_{2m}}
\newcommand{\CCC}{{\mathbf  C}}
\newcommand{\II}{{\mathrm{Im}}\,}
\newcommand{\RR}{{\mathrm{Re}}\,}
\newcommand{\eee}{{\mathrm  e}}
\newcommand{\LL}{L^2_\rho(\ren)}
\newcommand{\LLL}{L^2_{\rho^*}(\ren)}
\renewcommand{\a}{\alpha}
\renewcommand{\b}{\beta}
\newcommand{\g}{\gamma}
\newcommand{\G}{\Gamma}
\renewcommand{\d}{\delta}
\newcommand{\D}{\Delta}
\newcommand{\e}{\varepsilon}
\newcommand{\var}{\varphi}
\newcommand{\lll}{\l}
\renewcommand{\l}{\lambda}
\renewcommand{\o}{\omega}
\renewcommand{\O}{\Omega}
\newcommand{\s}{\sigma}
\renewcommand{\t}{\tau}
\renewcommand{\th}{\theta}
\newcommand{\z}{\zeta}
\newcommand{\wx}{\widetilde x}
\newcommand{\wt}{\widetilde t}
\newcommand{\noi}{\noindent}
\newcommand{\uu}{{\bf u}}
\newcommand{\xx}{{\bf x}}
\newcommand{\yy}{{\bf y}}
\newcommand{\zz}{{\bf z}}
\newcommand{\aaa}{{\bf a}}
\newcommand{\cc}{{\bf c}}
\newcommand{\jj}{{\bf j}}
\newcommand{\ggg}{{\bf g}}
\newcommand{\UU}{{\bf U}}
\newcommand{\YY}{{\bf Y}}
\newcommand{\HH}{{\bf H}}
\newcommand{\GGG}{{\bf G}}
\newcommand{\VV}{{\bf V}}
\newcommand{\ww}{{\bf w}}
\newcommand{\vv}{{\bf v}}
\newcommand{\hh}{{\bf h}}
\newcommand{\di}{{\rm div}\,}
\newcommand{\ii}{{\rm i}\,}
\newcommand{\inA}{\quad \mbox{in} \quad \ren \times \re_+}
\newcommand{\inB}{\quad \mbox{in} \quad}
\newcommand{\inC}{\quad \mbox{in} \quad \re \times \re_+}
\newcommand{\inD}{\quad \mbox{in} \quad \re}
\newcommand{\forA}{\quad \mbox{for} \quad}
\newcommand{\whereA}{,\quad \mbox{where} \quad}
\newcommand{\asA}{\quad \mbox{as} \quad}
\newcommand{\andA}{\quad \mbox{and} \quad}
\newcommand{\withA}{,\quad \mbox{with} \quad}
\newcommand{\orA}{,\quad \mbox{or} \quad}
\newcommand{\atA}{\quad \mbox{at} \quad}
\newcommand{\onA}{\quad \mbox{on} \quad}
\newcommand{\ef}{\eqref}
\newcommand{\mc}{\mathcal}
\newcommand{\mf}{\mathfrak}

\newcommand{\ssk}{\smallskip}
\newcommand{\LongA}{\quad \Longrightarrow \quad}
\def\com#1{\fbox{\parbox{6in}{\texttt{#1}}}}
\def\N{{\mathbb N}}
\def\A{{\cal A}}
\newcommand{\de}{\,d}
\newcommand{\eps}{\varepsilon}
\newcommand{\be}{\begin{equation}}
\newcommand{\ee}{\end{equation}}
\newcommand{\spt}{{\mbox spt}}
\newcommand{\ind}{{\mbox ind}}
\newcommand{\supp}{{\rm supp}\,}
\newcommand{\dip}{\displaystyle}
\newcommand{\prt}{\partial}
\renewcommand{\theequation}{\thesection.\arabic{equation}}
\renewcommand{\baselinestretch}{1.1}
\newcommand{\Dm}{(-\D)^m}

\title
{\bf Towards optimal regularity for the\\ fourth-order thin film
equation in $\re^N$:\\ Graveleau-type focusing self-similarity}


\author{P.~\'Alvarez-Caudevilla, J.D.~Evans, and V.A.~Galaktionov}

\address{Universidad Carlos III de Madrid,
Av. Universidad 30, 28911-Legan\'es, Spain -- Work phone number:
+34-916249099} \email{pacaudev@math.uc3m.es}

\address{Department of Mathematical Sciences, University of Bath,
 Bath BA2 7AY, UK -- Work phone number: +44 (0)1225386994 }
\email{masjde@bath.ac.uk}

\address{Department of Mathematical Sciences, University of Bath,
 Bath BA2 7AY, UK -- Work phone number: +44 (0)1225826988}
\email{masvg@bath.ac.uk}

\keywords{Thin film  equation in $\ren$, optimal regularity,
Graveleau-type focusing similarity solutions, numerics}

\thanks{This works has been partially supported by the Ministry of Economy and Competitiveness of
Spain under research project MTM2012-33258.}

 \subjclass{35G20, 35K65, 35K35, 37K50}

\date{\today}


\begin{abstract}

An approach to some ``optimal" (more precisely, non-improvable)
regularity of solutions of the \emph{thin film equation}
 $$
\tex{ u_{t} = -\nabla \cdot(|u|^{n} \nabla \D u)
 \quad \mbox{in} \quad \ren \times \re_+\,,
 \quad u(x,0)=u_0(x) \inB \re^N,}
 $$
 where $n \in (0,2)$ is a fixed exponent,
 with smooth compactly supported initial data $u_0(x)$,
 in dimensions $N \ge 2$  is discussed. Namely,
 a  precise exponent for the H\"{o}lder continuity with respect to the spatial radial variable $|x|$
is obtained by construction of a Graveleau-type focusing
self-similar solution. As a consequence, optimal regularity of the
  gradient $\nabla u$
 in certain $L^p$ spaces, as well as a H\"{o}lder continuity property of solutions with respect to
 $x$ and $t$,
 are derived, which cannot be obtained by classic standard
 methods of integral identities-inequalities. Several profiles for the solutions in the cases $n=0$ and $n>0$ are also plotted.

In general, we claim that, even for arbitrarily small $n>0$ and
positive analytic initial data $u_0(x)$, the solutions $u(x,t)$
cannot be better than  $C_x^{2-\e}$-smooth, where $\e(n)=O(n)$ as
$n \to 0$.

 \end{abstract}

\maketitle

\section{Introduction: TFE--4 and known related regularity results}
 \label{S1}


\noindent This paper is devoted to an ``optimal regularity"
analysis for a fourth-order quasilinear degenerate evolution
equation of the parabolic type, called the {\em thin film
equation} (TFE--4), with a fixed exponent $n>0$,
\begin{equation}
\label{i1}
    u_{t} = -\nabla \cdot(|u|^{n} \nabla \D u)
 \quad \mbox{in} \quad \ren \times \re_+\,,
 \quad u(x,0)=u_0(x) \inB \re^N,
\end{equation}
where $ n \in (0,2)$ is a fixed exponent,
and bounded, sufficiently smooth, and compactly supported initial
data $u_0$ with an arbitrary dimension
$N \ge 2$.

Moreover, we assume $u_0=u_0(|x|)$ to be radially symmetric, so
solutions $u=u(|x|,t)$ do the same. These initial conditions could
be relaxed (for example $u_0\in L^1\cap L^\infty$, etc.). However,
for simplicity, we have chosen those initial conditions, since the
focusing phenomenon to be studied here exists in any reasonable
functional class of initial data, including  radial ones.

The main result obtained here is an approach to the so-called
optimal regularity of solutions  of the TFE--4 for the
$N$-dimensional case with $N\geq 2$, ascertaining an ``optimal"
H\"{o}lder continuity exponent in $x$ and $t$.

For the one-dimensional
case, similar results were obtained in \cite{BF1}, establishing the
 H\"{o}lder continuity in $C^{0, \frac 12}$ with respect  to the variable $x$ and in $C^{0, \frac 1 8}$ with respect to $t$.
 However, in the $N$-dimensional case, the optimal regularity has been unsolved besides the
 interest of the specialised mathematical community;
 see \cite{BHRR,DGG, Grun95, Grun04}. In many of those works it is also assumed
 that the solutions are non-negative
  for the Cauchy problem (CP). However, as we recently proved in the CP settings, the
  solutions for the TFE--4 \eqref{i1} are oscillatory and sing-changing for all not that
  large values of $n>0$;
 see   \cite{PV3, EGK2} for recent results and references therein.

One of the difficulties concerning optimal regularity for the
equation \ef{i1} is due to the fact that this equation \eqref{i1} is {\em not fully divergent}, so the
whole ``global" regularity theory for the TFE--4 can use just a
couple of well-known integral identities/inequalities. Indeed, as discussed in \cite{PV3} an integral identities argument,
as that performed by Bernis--Friedman in 1D, fails to show the H\"{o}lder continuity in higher dimensions; see \cite{PV3} for details and reasons for these issues.
Consequently, we claim
that, for such ``partially" (and/or fully non-) divergent
equations, a different approach should be put in charge.


 \subsection{Methodology: Graveleau-type focusing similarity solutions}


 Here, working on the radially symmetric problem of TFE--4,
  such that
  \[u(x,t)=u_*(r,t),\quad \hbox{with} \quad r=|x|,\]
  we base our analysis on the
 \emph{focusing argument} performed, firstly,  by Graveleau \cite{Grav72} in a formal sense, and later justified  in
  Aronson--Graveleau \cite{AG} for a class of non-negative solutions (via the Maximum Principle, MP)
  of the {\em classic porous medium
  equation} (PME--2):
 \begin{equation}
\label{pormed}
    v_{t} = \D (v^m)
 \quad \mbox{in} \quad \ren \times (-\iy,0),\quad \hbox{with an exponent $m>1$.}
\end{equation}
Actually, in comparison with \ef{i1}, $m=1+n>1$. In particular, the authors in \cite{AG} obtained,
after constructing a one-parameter family of self-similar
solutions, that,
$$\hbox{when}\quad  N\geq 2, \quad \hbox{there exists}\quad 0<\mu=\mu(m,N)<1,$$
and a radial self-similar solution $v_*(r,t)$ to the focusing
problem and with the focusing time $t=0^-$, such that\footnote{We do
not put further details how to get \ef{Hol1}, since we will explain
this shortly for the TFE--4.}
\be
\label{Hol1}
 v_*(r,0^-)=C_* r^\mu,\quad \hbox{where $C_*$ is an arbitrary
constant.}
 \ee
 This finished a long discussion concerning optimal regularity for
 the PME--2 after seminar results by Caffarelli--Friedman
 \cite{Caff79} (see {\tt MathSciNet} for many further extensions)
 proving, earlier, H\"older continuity of $v(x,t)$ in a very
 general setting.

 The above results of Aronson and Graveleau
showed that such a H\"older continuity is, indeed, optimal, and
even the corresponding H\"older exponent $\mu=\mu(m,N)$ from
\ef{Hol1} can be evaluated, at least, numerically for any $m>1$
and $N \ge 2$. Indeed, the blow-up singularity \ef{Hol1} shows,
actually, not an {\em optimal} regularity for such PME--2 flows,
but the {\em non-improvable} one, in the sense that, for the
H\"older exponent of solutions of \ef{pormed},
 \be
 \label{Opt1}
 \mbox{a non-improvable regularity} \quad \Rightarrow \quad
 \mbox{H\"older exponent $ \not >$ $\mu$}.
 \ee
  Fortunately,
as we have mentioned above, Caffarelli--Friedman earlier results
of the 1980s had been already available at that time, so,
together with \ef{Opt1}, those proved that a proper H\"older
regularity of any non-negative solutions of the PME--2 is optimal.

These Aronson--Graveleau focusing (singularity--blow-up) ideas for
the PME--2 \ef{pormed} and other second-order quasilinear
parabolic equations (e.g. for  the $p$-Laplacian one) got further
development and extensions in a number of papers devoted to
various important aspects of porous medium flows; see references
and results in \cite{Ang01, Ar03}.

   More precisely talking about such Graveleau-type similarity approach, from the point of view
   of applications, in the focusing
problem, one assumes that, initially, there is a non-empty compact
set $K$ in the complement of the ${\rm supp}\,v_0$, where $v$
vanishes. In other words, there is a hole in the support of the
initial value $v(x,0)=v_0(x) \ge 0$, and, in finite time $T$, this
hole disappears, making the solution $v$ to become positive along
the boundary of $K$ and eventually at all points of $K$.
Basically,  due to the finite propagation property that these
equations possess (the porous medium equation and the thin film
equation). Thus, as the flow evolves the liquid enters $K$ and
eventually reaches all points of $K$ at the instant $T$. We are
then interested in how the solution for the TFE--4 \eqref{i1}
behaves near the focusing time $T$. We again suppose (as in the
Graveleau's argument for PME--2) the focusing time $t=0^-$.

However, in our case, with the TFE--4 \eqref{i1} and $n>0$, but
 close to zero, we need to take into consideration the existence of
oscillatory solutions of changing sign.

Furthermore, for the TFE--4 \ef{i1}, there are not still any
general regularity results even in the radial setting in $\ren$.
Consequently, we will follow the lines and the logic  of \ef{Opt1},
i.e. we will estimate a certain {\em non-improvable} H\"older
exponent for radial (and, hence, all other) solutions.

To this end, we use those previous singularity (blow-up) ideas
from the PME--2 \eqref{pormed} to establish some properties for
the self-similar solutions of the radially symmetric problem of
the TFE--4 \eqref{i1} and, eventually, some optimal regularity
information about its more general solutions. Indeed, as we will
see through the analysis performed in this paper, to obtain such
an ``optimal" (non-improvable) regularity for the TFE--4
\eqref{i1},  it seems that we must solve a nonlinear focusing
problem for that TFE--4 that will be derived from the associated
self-similar equation.

To be precise, we first apply the focusing ideas performed by Aronson-Graveleau \cite{AG} for the PME--2 \eqref{pormed} to the TFE--4 \eqref{i1}.
To do so, we will work on the space of radially self-similar solutions with $r=|x|>0$ of the form
\be
\label{sssol}
 \tex{
u_*^{\pm}(r,t)= (\pm t)^{\a} f(y), \quad y=\frac{r}{(\pm t)^\b}
\,\,\, \mbox{for}
  \,\,\, \pm t > 0 \whereA \b= \frac {1+\a
 n}4>0.
 }
\ee
Thus, these solutions will solve an associated self-similar equation, or \emph{nonlinear eigenvalue problem},
$$
 \tex{
- \frac{1}{y^{N-1}} \big[ y^{N-1}|f|^n
\big(\frac{1}{y^{N-1}}(y^{N-1} f')'
 \big)' \big]' \pm \b y f' \mp \alpha f = 0,
 }
 $$
 with $-$ or $+$ depending on if we are analysing the focusing phenomena or the defocusing phenomena. In other words, before
 the focusing time or after (in our situation $t=0$).

 As usual, the unknown {\em a priori}
  $\a$ represents the nonlinear eigenvalues supporting the reason to call that equation a \emph{nonlinear eigenvalue problem}.

 Crucial in our analysis will be to impose some radiation-type conditions, or some minimal and maximal growth
  at infinity, i.e. as $y\to +\infty$. This
 actually allows us, among other things, to determine the existence of a discrete family of nonlinear eigenvalues denoted by $\{\a_k\}$.

 Moreover, this behaviour at infinity leads us to have
  the solutions of the self-similar TFE--4 bounded at infinity by functions of the form
$$
 \tex{
f(y)=C  y^{\mu} (1+o(1)) \quad
C \in \re \quad \hbox{and}\quad \mu=\mu(\a,n):=\frac{4\a}{1+\a n} >0.
 }
$$ Thus, we ascertain a limit at the focusing point
 $$
  \tex{
  u_*(r,t) \to
C r^{\mu} \asA t \to 0^-, \quad \mbox{with} \quad
  \mu= \frac \a \b,
 }
  $$
  uniformly on compact subsets (see details in Section \ref{S3}).
  However, this does not provide yet information about the influence of the dimension on the behaviour of the solutions at infinity.

  To get to that point, we need to find those nonlinear eigenvalues\footnote{Since
 we obtain those nonlinear eigenvalues as a perturbation from a family of linear eigenvalues when the parameter $n$
  is  close to zero, we will
 use indistinctively this notation for the family of linear or nonlinear eigenvalues. To be precise,
  the nonlinear eigenvalues correspond to the family $\{\a_k(n)\}$ and
 the linear eigenvalues to the family $\{\a_k(0)\}$.} $\{\a_k\}$ associated with a discrete family
  of nonlinear eigenfunctions $\{f_k\}$. Here, we have done this analysis via a homotopy deformation as $n\to 0^+$ (Sections \ref{S.lin} and \ref{S5})
  to the linear problem with $n=0$,
  \be
  \label{linbh}
  u_t=- \D^2 u \inB \ren \times \re_-,
  \ee
  and with those patterns occurring at $n = 0$ acting as branching points of nonlinear eigenfunctions of
  the Cauchy problem \eqref{i1} when $n$ is  close to zero.
   In particular, according to our analysis, when $n=0$, it follows that
 $$
  \tex{
  \a_k(0)=\a_k=\frac{k}{2},\quad k=1,2,3,\cdots.
 }
  $$
 As an observation, those values of the parameter $\a$ when $n=0$ might be written as
  a perturbation of the eigenvalues for the fourth-order operator
 $$
 \tex{
 {\bf B}^* f \equiv - \frac{1}{y^{N-1}} \big[ y^{N-1} \big(\frac{1}{y^{N-1}}(y^{N-1} f')'
 \big)' \big]' -  \frac 14\,  y f'=\l f,
 }
 $$
analysed in full detail in \cite{EGKP}; see further comments below.

Hence, we perform a homotopy deformation from the TFE--4
\eqref{i1} to the parabolic bi-harmonic equation \eqref{linbh},
for which we can apply again a similar logic in ascertaining
Graveleau-type ``focusing solutions". Indeed, we get a minimal and
maximal growth for the radial self-similar solutions of this
linear problem \eqref{linbh}. Moreover, since we know that there
exists a discrete family of eigenvalues for the corresponding
self-similar equation associated with \eqref{linbh}, using this
branching/homotopy argument we obtain the desired family of values
for the parameter $\a$ from which we will have a  family of radial
self-similar solutions emanating at $n=0$ from the self-similar
solutions of the linear problem \eqref{linbh}. Note also that  the
family of eigenfunctions for the linear problem \eqref{linbh} is a
complete set of generalised Hermite polynomials with finite
oscillatory properties \cite{EGKP}.

In this paper, we also  support this analysis numerically,
performing a shooting procedure from $y=0$ to $y=+\infty$ that
gives us those linear eigenvalues $\a_k$ and the profile of the
corresponding eigenfunctions.

Additionally, we show some numerical analysis that provided us,
with very difficult to ascertain, profiles of the nonlinear
eigenfunctions $\{f_k\}$ with $n>0$ and sufficiently close to
zero.

This numerical analysis supports the conjecture that by  homotopy
continuity the properties for the profile when $n=0$ remain valid
for $n>0$ and sufficiently close to zero. As seen above, we have
fixed the interval $n\in (0,2)$, however, the numerics suggest
that those properties are lost from $n=0.8$ or $0.9$. In fact, if
$n\geq 2$ we just have nonnegative solutions as shown in
\cite{BF1} so that, this focusing argument is not possible. Thus,
we restrict this analysis to $n\in (0,2)$, or more precisely, of
$n>0$ sufficiently close to zero.

 Finally, we obtain such an ``optimal" (non-improvable) regularity
for the TFE--4 \eqref{i1} proving that there is H\"{o}lder
continuity with an specific coefficient $$\mu_k
=\frac{4\a_k}{1+n\a_k}\quad \hbox{and}\quad \n u(x,t) \in
L^{p}_{\rm loc}(\ren) \quad \hbox{with}\quad p <
p^*(n,N)=\frac{N}{1-\mu_k},$$ having such regularisation depending
on the dimension $N$, as we claimed and were looking for.
Consequently, if $n>0$ small,  we find that
\be
\label{reguc3} u_*(r,t) \in C_r^{2-\e(n)},\quad \hbox{with} \quad
\e(n)=n\mu+O(n^2). \ee It turns out that, thanks to our branching
analysis, there holds $\mu>0$, so that, for $n>0$, the regularity
condition \eqref{reguc3} becomes obviously worse. Note that, if
$\mu<0$, we have a better regularity than for $n=0$ which is
impossible.

In general, we obtain (and the numerics suggest that as well) that
there exist several focusing singularities in the radial geometry
of the type $C^{2k-\e}$ with $k=1,2,3,\cdots$. In particular, for
the values of $\a_k$ obtained here, it follows that for $k=1$,
$\a_1=\frac 1 2$ is the minimal and crucial one, since it seems to
be the only degenerate case, i.e.  $f_1(0)=0$. All others satisfy
$$f_{k}(0)>0, \quad K=1,2,3,\cdots,\quad \hbox{i.e.  the TFE--4
for $n>0$ is not degenerate}
 $$
 initially, but, indeed, will be eventually, at the focusing time $t=T^-(=0^-)$. Therefore, we can conclude that slight
changes in the parameter $\a$ destroy the regularity, supporting
the fact that these results are valid only when $n$ is
sufficiently close to zero. Indeed, via branching analysis we find
that $$\a_k(n)=\a_k+\mu_k n+O(n^2),$$ so that we have the
regularity condition $$u_*(r,t) \in C^{2k-\mu_k n+O(n^2)}.$$ Note
that, even for $\a_2=1$, the analytic positive solution becomes,
at $t=T^-$, $C^{4-\e}$, i.e. not classical solutions in $C^4$.
Hence, for the minimal $\a_1=\frac 1 2$ it is $C^{2-\e}$, which is
much worse.

\vspace{0.2cm}

\noindent{\bf Remark.} As mentioned above, we also observe that,
for the PME--2 \eqref{pormed}, one needs to have such a non-empty
compact hole to apply the Maximum Principle. However, we believe
that, for the TFE--4 \eqref{i1} and the value $\a_1=\frac 1 2$,
that provides us with the minimal regularity $C^{2-\e(n)}$, there
is not such a hole. We do not have a rigorous justification of it
but the numerics presented in this paper suggest it that way.

\vspace{0.2cm}

In relation to the H\"{o}lder regularity with respect to the temporal variable $t$, 
assuming the radially self-similar solutions of the form \eqref{sssol} we find that
$$|u(0,t)-u(0,0)| = (-t)^\a f(0).$$
Hence,the H\"{o}lder's exponent for the variable $t$ when it is very close to $t=0$ cannot be bigger than 
$$\a_1(n)=1/2 + O(n).$$
Therefore, we provide with an estimation from above for the H\"{o}lder continuity with respect to the $t$ variable. Moreover,
this estimation improves the H\"{o}lder's exponent obtained by Bernis--Friedman \cite{BF1} showing that theirs was actually not optimal since it was $1/8$.

\subsection{TFE-4 problem settings}

This setting is well-known nowadays (though many things were not
fully proved in a general case), so, below, we omit many details;
see surveys  \cite{Grun95, Grun04}. Note that, for some values of
$n>0$ (not that large), focusing similarity solutions to be
constructed do not exhibit finite interfaces, so the resulting
optimal regularity results are true for any FBP and/or Cauchy
problem settings. Principal differences between the CP and the
standard FBP settings for the TFE--4 \ef{i1} are explained in
\cite{EGK1, EGK2}.

We recall that the
solutions are assumed to satisfy the following zero contact angle boundary
conditions:
\begin{equation}
\label{i5}
    \left\{\begin{array}{ll}  \tex{ u=0,} & \tex{ \hbox{zero-height,} }  \\
    \tex{ \nabla u=0,} & \tex{ \hbox{zero contact angle,} } \\
    \tex{ -{\bf n} \cdot \nabla (\left|u\right|^{n}  \D u)=0,}  &
    \tex{ \hbox{conservation of mass (zero-flux)} } \end{array} \right.
\end{equation}
at the singularity surface (interface) $\Gamma_0[u]$, which is the lateral boundary
of
\begin{equation*}
    \tex{ \hbox{supp} \;u \subset \re^{N} \times \re_+,\quad N \geq 1\,,}
\end{equation*}
where ${\bf n}$ stands for the unit outward normal to
$\Gamma_0[u]$, which is assumed to be sufficiently smooth (the treatment of such hypotheses is
not any goal of this paper).
For  smooth interfaces, the condition on the flux can be read as
\begin{equation*}
    \lim_{\hbox{dist}(x,\Gamma_0[u])\downarrow 0}
    -{\bf n} \cdot \nabla (|u|^{n}  \D u)=0.
\end{equation*}
Next, we denote by
\begin{equation}
\label{mass}
 \tex{
    M(t):=\int  u(x,t) \, {\mathrm d}x
    }
\end{equation}
the mass of the solution, where integration is performed over the support. Then,
differentiating $M(t)$ with respect to $t$ and applying the
divergence theorem, we have that
\begin{equation*}
 \tex{
  J(t):=  \frac{{\mathrm d}M}{{\mathrm d}t}= -
  \int\limits_{\Gamma_0\cap\{t\}}{\bf n} \cdot \nabla
     (|u|^{n}  \D u )\, .
     }
\end{equation*}
The mass is conserved if $ J(t) \equiv 0$,
which is satisfied by the flux condition in \eqref{i5}.

\section{Radial self-similar solutions: focusing and defocusing cases}
 \label{S2}


\noindent Now, we construct the operators and specific solutions in order to apply
the ideas performed by Aronson--Graveleau \cite{AG} for the PME--2 \eqref{pormed} to the TFE--4 \eqref{i1}.

Thus, thanks to the scaling invariant property of these nonlinear
parabolic equations, we now construct radially self-similar
solutions, i.e. in terms of $r=|x|>0$, with a still unknown value
of the parameter $\a>0$ (clearly, it must be positive, as shown
below)
\be
  \label{upm}
   \tex{
 u_*^{\pm}(r,t)= (\pm t)^{\a} f(y), \quad y=\frac{r}{(\pm t)^\b} \,\,\, \mbox{for}
  \,\,\, \pm t > 0 \whereA \b= \frac {1+\a
 n}4>0.
  }
  \ee
Here, by the time-translation, we ascribe the blow-up  or
\emph{focusing} time to $T = 0^-$.
 We then simultaneously consider two cases:
\begin{enumerate}
\item[(i)] {\bf Focusing} (i.e., Graveleau-type) similarity solutions, which
 play a key role, corresponding to $(-t)$ in \ef{upm} and the
 singular blow-up limit as $t \to 0^-$, and

\item[(ii)] {\bf Defocusing} similarity solutions, with $(+t)$ in
 \ef{upm}, playing a secondary role as extensions of the previous
 ones for $t>0$, i.e.,  corresponding to the not-that-singular (or,
 at least, less singular) limit $t \to 0^+$. Actually, these
 defocusing solutions are well-posed solutions of the Cauchy
 problem for the TFE--4 for $t>0$ with initial data $u_*(r,0^-)$
 obtained from the previous blow-up limit as $t \to 0^-$.
\end{enumerate}

 \ssk

Substituting \eqref{upm} into \eqref{i1} we arrive at the
similarity profiles $f(y)$ satisfying the following {\em
  nonlinear eigenvalue problems}:
   \be
   \label{eigpm}
   \tex{
 \BB^\pm_{n} (\a,f) \equiv - \n_y \cdot (|f|^n \n_y \D_y f) \pm \b y \cdot \n_y f
 \mp
 \a f=0 \,\,\,\mbox{for} \,\,\, y>0; \quad f'(0)=f'''(0)=0,
 }
 \ee
 where $\nabla_y$ and $\D_y$ stand for the radial gradient and the radial Laplacian.
 In \ef{eigpm}, we present two symmetry conditions at the origin
 (which should be modified if $f(y) \equiv 0$ near the origin,
 i.e. it  contains a ``zero hole" nearby). As we will see with the numerical analysis
 performed in Section \ref{S.lin} of this paper, we can choose other conditions depending on if either $f=0$ or $f\neq 0$.

 To complete these
 nonlinear eigenvalue settings one needs extra ``radiation-type"
 (or growth-type) conditions at infinity to be introduced next.

  Indeed, we actually find
 radially similarity profiles $f$ depending on the single variable $y$. The operator of the equation \eqref{eigpm} is then
\be
 \label{rad22}
 \tex{{\bf A}^\pm_{n} (\a,f) \equiv- \frac{1}{y^{N-1}} \big[ y^{N-1}|f|^n \big(\frac{1}{y^{N-1}}(y^{N-1} f')'
 \big)' \big]' \pm \b y f' \mp \alpha f = 0,}
 \ee
 denoting the radial nonlinear operator as
 \be
 \label{rad23}
  \tex{
 {\bf R}_{n}(\a,f)= \frac{1}{y^{N-1}} \big[ y^{N-1}|f|^n \big(\frac{1}{y^{N-1}}(y^{N-1} f')'\big)' \big]'.
 }
 \ee

 Thus, here, $\a >0$ is  a parameter, which, in fact, stands in both cases for admitted real
  {\em nonlinear eigenvalues} to be determined, in the focusing Graveleau-type case, via the
  solvability of the corresponding {\em nonlinear eigenvalue
  problem}, accomplished with some special (radiation-type)
  conditions at infinity, to be introduced shortly.

\begin{remark}
{\rm In general, with respect to the similarity profiles described above we note that there are two main types of self-similar
solutions. For solutions of the first kind the similarity variable
$y$ can be determined \emph{a priori} from dimensional
considerations and conservation laws, such as the conservation of
mass \eqref{mass} or momentum.

For solutions of the second kind the exponent $\b$ (and by relations the exponent $\a$) in the similarity variable
must be obtained along with the solution by solving a nonlinear eigenvalue problem of the form \eqref{eigpm}.

The first published examples of self-similar solutions of second
kind are due to G.~Guderley in 1942 \cite{G42} studying imploding
shock waves, although the term was introduced by Ya. B.~Zel'dovich
in 1956 \cite{Zel56}.  Further examples might be found in the theory of the collapse of bubbles in
compressible fluids or in works on gas motion under an impulsive load; see Barenblatt \cite{B}
 for an extensive work on this matter.
}
\end{remark}

\section{Minimal growth at infinity  (a ``nonlinear radiation condition") for Graveleau-type profiles}
 \label{S3}


\noindent 
Here, we consider the blow-up problem $(\ref{rad22})_-$, i.e.
with the lower signs. One concludes that, in order to get a
possible discrete set of eigenvalues $\a>0$, some extra conditions
denoted by nonlinear radiation condition on the behaviour (in
particular, growth) of $f(y)$ as $y \to + \iy$ must be imposed.

Obviously, such a ``radiation-type"
condition (we use a standard term from dispersion theory) just
follows from analysing all possible types of such a behaviour,
which can be admitted by the ODE $(\ref{rad22})_-$, which is not
that a difficult problem.
There are two cases:

\begin{enumerate}
\item[(I)] {\em More difficult: there is a ``zero hole" for $f(y)$ near the origin.}
 Finite interfaces for the TFE--4 are well known in the case of
 the standard FBP setting and also in the Cauchy one (see
 \cite{EGK1,EGK2} and references therein).
 Then  we need to look for
profiles $f(y)$ which vanish at finite $y=y_0>0$ and describe
asymptotics of the general solution satisfying zero contact angle
and zero flux conditions as $y\to y_0^+$
  $$
  f(y)\to 0,\quad
f'(y)\to 0,\quad -|f|^n f'''(y)\to 0.
 $$
This case causes a difficult problem, since the oscillatory
behaviour of $f(y)$ close to the interface is quite tricky for the
CP \cite{EGK1}, while, for the FBP, it is much better understood.

\item[(II)] {\em More standard and easier: $f(0) \ne 0$.} This case also
includes the border possibility $y_0=0$, i.e.   when $f(0)=0$,
but $f \not \equiv 0$ in any arbitrarily small neighbourhood of
$y=0$.

\end{enumerate}


 \subsection{Minimal and maximal growth at infinity}

This is key for our regularity analysis.
  Our radial ODE
$(\ref{rad22})_-$ admits two kinds of asymptotic behaviour at
infinity, i.e.  when $y_0 \to +\infty$. For the operator ${\bf
A}^-_{n} (\a,f)$, we state the following result:

\begin{proposition}
\label{theneg} For any $\a>0$, the  ODE $\eqref{rad22}_-$, with
the operator ${\bf A}^-_{n} (\a,f)$  possesses:

{\rm (i)} Solutions $f(y)$ with a minimal growth
  \be
  \label{gg1}
  \tex{
f(y)=C  y^{\mu} (1+o(1))\quad \hbox{as}\quad y\to +\infty,} \quad
C \in \re,
 \ee
  where
\be
\label{param}
 \tex{
 \mu=\mu(\a,n)=\frac{4\a}{1+\a n} >0.
 }
   \ee

{\rm (ii)} Moreover, there exist
 solutions of
 $\eqref{rad22}_-$ with a maximal growth
 \be
 \label{gg2}
 \tex{
  f(y) \sim
y^{\mu_0},\quad \hbox{as}\quad y\to +\infty, \quad \mbox{where}
\quad \mu_0= \frac 4 n.
 } \ee
\end{proposition}

\vspace{0.2cm}

\noindent{\bf Remark.}
\begin{itemize}
\item It is important to mention that the expression \ef{gg2} for
solutions of the maximal growth does not include the ``additional or extra" corresponding
{\em oscillatory component} $\var(s)$, with $s= \ln y$, and shows
just the growth behaviour of its ``envelope" with algebraic
growth. Such an oscillatory maximal behaviour will be introduced
and studied later.

\item Note also that, as $n
\to 0^+$, we have that
 \be
 \label{gg81} \tex{
  \frac 4n \to +\iy,
  }
  \ee
  which corresponds to an exponential oscillatory growth for the
  linear problem occurring for $n=0$; see next section.

\item The difference between \ef{gg1} and \ef{gg2}, which implies such
terms as minimal and maximal growth at infinity, is obvious:
 \be
 \label{gg3}
 \tex{
  \mbox{for any $\a>0$}, \quad \frac 4n >  \frac \a \b = \frac {4
  \a}{1+\a n}.
  }
  \ee
  \end{itemize}

  \vspace{0.2cm}

 {\em Proof.}
 This follows from a balancing of linear and
 nonlinear operators in this ODE, though a rigorous justification
 is rather involved and technical. A formal derivation is surely
 standard and easy:

\begin{enumerate}
\item[(i)] In this first case, we assume a linear asymptotic as $y\to +\infty$
(assuming simple radial behaviour $f \sim y^{\mu}$)
\be
  \label{kk3}
  \tex{
  \frac{1+\a n}4 \, y  f' - \a f+...=0
  \LongA f(y) \sim y^\mu \whereA \mu= \frac {4\a}{1+\a n}>0.
   }
   \ee
   Since this behaviour is asymptotically linear, eventually, we
   get an arbitrary constant $C \not = 0$ in \ef{gg1}.
 A full justification of existence of such orbits is
 straightforward, since the nonlinear term in the ODE is then
 negligible. So that in the equivalent integral equation
 though a singular term, it produces a negligible perturbation.

 \item[(ii)]   On the other hand, the solutions are also bounded by a maximal growth, which in this case comes from non-linear asymptotics
 that balance all  three  operators
    \be
 \label{kk2}
  \tex{
  {\bf R}_{n}(f)+ \frac{1+\a n}4\, y f' - \a f=0
  \LongA
  f(y) \sim y^{\frac 4n} \asA y \to +\iy,
  }
 \ee
where we again indicate the envelope behaviour of this oscillatory bundle.
A justification here, via Banach's contraction principle is even
easier, but rather technical, so we omit details.
 \qed
 \end{enumerate}

\ssk

Overall this allows us to formulate such a condition at infinity,
which now takes a clear ``minimal nature" such that solutions $f$
are now bounded at infinity by a function
\be
\label{mingro} \tex{
  f(y)=C   y^{\mu} (1+o(1)),\quad
\hbox{with}\quad \mu= \frac {4\a}{1+\a n}>0.} \ee
Obviously, thus we just need a global solution $f(y)$ of our ODE
$(\ref{rad22})_-$ in $\re_+$, satisfying the minimal growth
\ef{gg1}. Indeed, for such profiles, there exists a finite limit
at the focusing (blow-up) point:
 \be
 \label{gg5}
  \tex{
  u_*(r,t) \to C r^{\mu} \asA t \to 0^-, \quad \mbox{with} \quad
  \mu= \frac \a \b,
  }
  \ee
 uniformly on compact intervals in $r=|x| \ge 0$. One can see
 that, for any maximal profile as in \ef{gg2}, the limit as in
 \ef{gg5} is infinite, so that such similarity solutions do not
 leave a finite trace as $t \to 0^-$.

However, the above proposition leaves aside the principle question
on the dimensions of the corresponding minimal and maximal bundles
as $y \to +\iy$, which is not a straightforward problem. Note that
the latter is actually supposed to determine the strategy of a
well posed shooting of possible solutions of the above focusing
problem.

Indeed, the main problem is how to find those admitted values of
nonlinear eigenvalues $\{\a_k\}_{k \ge 1}$ (possibly and
hopefully, a discrete set), for which $f=f_k(y)$ exist producing
finite limits as in \ef{gg5}.

To clarify those issues, we consider the much simpler linear case when
$n=0$ and, subsequently, pass to the limit as $n\to 0^+$ in \eqref{rad22}. This analysis will provide us,
eventually, with some qualitative information about the solutions, at least when $n$ is very close to zero.


\section{The linear problem: discretisation of $\a$ for $n=0$}
 \label{S.lin}

 For $n=0$, the TFE--4 \ef{i1} becomes the classic {\em
 bi-harmonic equation}
  \be
  \label{bi.1}
  u_t=- \D^2 u \inB \ren \times \re_-.
   \ee
 Of course, solutions of \ef{bi.1} are analytic in both $x$ and
 $t$, so any focusing for it makes no sense.
 However, we will show that (following a similar philosophy to the one above) Graveleau-type ``focusing solutions"
 for \ef{bi.1} are rather helpful to predict some properties of
 true blow-up self-similar solutions of \ef{i1}, at least, for
 small $n>0$.

\subsection{Maximal and minimal bundles}

 Thus, first we consider the same ``focusing" solutions $(\ref{upm})_-$ for
\ef{bi.1}, that take a simpler form
  \be
  \label{b1}
   \tex{
 u_*(r,t)=(-t)^\a f(y), \quad y = \frac   r{(-t)^{1/4}} \quad  \big(\b=
 \frac 14  \big).
  }
 \ee
  Then, the corresponding linear radial ODE \ef{rad22} takes also a
  simpler form
\be
 \label{b2}
 \tex{{\bf A}^-_{0,y} (\a,f) \equiv- \frac{1}{y^{N-1}} \big[ y^{N-1} \big(\frac{1}{y^{N-1}}(y^{N-1} f')'
 \big)' \big]' -  \frac 14\,  y f' + \alpha f = 0.}
 \ee
We first calculate the solutions of \ef{b2} with a maximal
behaviour. Those are exponentially growing solutions of the form
 \be
 \label{b3}
  \tex{
  f(y) \sim \eee^{a y^\g} \asA y \to + \iy \LongA a^3=- \frac 14 \,\big(
  \frac 34 \big)^3.
   }
   \ee
 This characteristic equation gives two roots with ${\rm Re}(\cdot)>0$:
  \be
  \label{b4}
   \tex{
  a_{1,2}=\frac 34\,4^{-\frac 13} \big[ \frac 12 \pm {\rm i} \, \frac{\sqrt 3}2
  \big] \equiv c_0 \pm {\rm i} \, c_1.
 }
 \ee
 and one negative root (this actually goes to the bundle of minimal solutions;
 see below)
  \be
  \label{b5}
   \tex{
 a_3= -\frac 34 \, 4^{-\frac 13}<0.
 }
  \ee
 Hence, the bundle of maximal solutions is oscillatory as $y \to +
 \iy$ (including the multiplicative algebraic factor),
  \be
  \label{b51}
  \tex{
  f(y) \sim y^{-\frac{2}{3}(N+2\alpha)} \eee^{c_0 y^{4/3}}[C_1 \cos (c_1 y^{\frac 4 3}) + C_2
  \sin(c_1 y^{\frac 43})], \quad C_{1,2} \in \re.
 }
 \ee
 In fact, obviously, this bundle is 3D since it includes the
 1D sub-bundle of exponentially decaying solutions with the
 exponent \ef{b5}. However, in a shooting procedure performed at the end of this section, we
 intend to get rid of just two coefficients in \ef{b51}
 $$C_1=C_2=0.$$
 Indeed, this numerical evidence
 shows that those coefficients vanish at certain values of the parameter $\a$ that we will determine below specifically.
 Moreover, we provide an explanation in the next subsection.

On the other hand, the minimal behaviour \ef{gg1} now reads
 \be
 \label{b6}
  \tex{
  f(y)= C y^{4 \a}(1+o(1)), \quad C \in \re \quad \big( \hbox{with}\quad \frac \a \b= 4 \a \big).
   }
   \ee
The whole bundle of such minimal solutions is 2D. Besides the
parameter $C$ in \ef{b6}, it includes a 1D sub-bundle of
exponentially decaying solutions with the exponent \ef{b5}, so
that, overall, the minimal solutions compose a 2D family:
 \be
 \label{b7}
  \tex{
  f(y) \sim C y^{4 \a}(1+o(1))+ D \eee^{a_3 y^{4/3}} (1+o(1)),
  \quad C,D \in \re.
  }
  \ee
 Justification of both behaviours is, indeed, simpler for such a standard linear ODE problem.

Finally, passing to the limit as in \ef{gg5} then gives, for
minimal solutions, a finite ``focusing trace":
 \be
  \label{b8}
 u_*(r,t) \to C r^{4\a} \asA t \to 0^-.
  \ee
We explain why those trivial results are so important in
what follows. The conclusion from \ef{b8} is straightforward:
since $u_*(r,0^-)$ must be analytic then we have the following.

 \begin{proposition}
  \label{Pr.anal}
  For the above linear ``focusing eigenvalue problem" $\ef{b2}$,
   there exists not more than  a countable set of admissible
  eigenvalues $\{\a_k\}_{k \ge 1}$ given by
   \be
   \label{b9}
    \tex{
    \a_k= \frac k2, \quad k=1,2,3,...\,.
    }
    \ee
    \end{proposition}

     {\em Proof.}
     The function $r^{4\a}$ is analytic at
    $r=0$ only for values $\a_k$ from \ef{b9}, i.e. $4 \a$ must
    be a real even integer, say $2k$. Note also  that, in the
    present problem, all functions are analytic, so that,
    automatically, the set of roots $\{\a_k\}$ is discrete with a
    possible accumulation point at infinity only.
     \qed

\ssk

\noindent{\bf Remark.} Solving the {\em linear eigenvalue
problem}:
\[
 \tex{
 {\bf A}^-_{0,y} (\a_k,f_k) = 0 \inB \re, \quad f_k \in L^2_\rho(\re),
  }
  \]
  then, it seems that the nonlinear eigenvalue
   problem $\eqref{rad22}_{-}$ formally reduces to the classic linear eigenvalue problem \eqref{b2} at $n=0$,
   providing us with another reason to call  \eqref{rad22} a {\em nonlinear eigenvalue
problem}.

 Also, the values of the parameter $\a$ can be written as shifting from the eigenvalues of the eigenvalue problem
 \[
 \tex{
  {\bf B}^* f \equiv - \frac{1}{y^{N-1}} \big[ y^{N-1} \big(\frac{1}{y^{N-1}}(y^{N-1} f')'
 \big)' \big]' -  \frac 14\,  y f'=\l f,
 }
 \]
 analysed  in \cite{EGKP} and whose discrete spectrum takes  the form
 \[
 \tex{
 {\bf \s}({\bf B}^*)=\{\l_k=-\frac{k}{4}\,;\, k=0,1,\cdots\}.
 }
 \]
 Hence,
 \[
 \tex{
 \a_k=\l_k -\frac{k+1}{4},\quad \hbox{with}\quad k=1,2,3,\cdots,
 }
 \]
 having a countable family of eigenvalues for the problem \eqref{b2}.

 Moreover, note that ${\bf A}^-_{0,y} (\a,f)$ \eqref{b2} is a non-symmetric linear operator, which is bounded
from $H_{\rho}^4(\re)$ to $L_{\rho}^2(\re)$ with the exponential
weight given by \eqref{b51}, i.e., $$\rho(y)= {\mathrm e}^{-c_0
|y|^{4/3}}, \quad c_0>0 \,\,\,\mbox{small and defined by
\eqref{b4}}.$$


 \subsection{Well-posed shooting procedure}

 We now discuss a practical procedure to obtain the linear eigenvalues
 $\a_k$. We perform standard shooting from $y=0$ to
 $y=+\iy$ for the ODE \ef{b2} by posing four conditions at the
 origin:
 \begin{itemize}
\item Either
  \be
  \label{sh1}
  f(0)=1 \,\,\,\mbox{(normalisation)}, \quad f''(0)=0, \quad f'(0)=f'''(0)=0
  \,\,\,\mbox{(symmetry)};
   \ee
 \item or \be
  \label{sh2}
  f(0)=0, \,\,\, \quad f''(0)=1 \,\,\, \mbox{(normalisation)}, \quad f'(0)=f'''(0)=0
  \,\,\,\mbox{(symmetry)}.
   \ee
   \end{itemize}
 These two sets of conditions correspond to partitioning of the eigenfunctions into two subsets with
 the corresponding properties at the origin of
 $$\hbox{either}\quad f(0)=1,\quad f''(0)=0,\quad \hbox{or}\quad f(0)=0, \quad f''(0)=1,$$
 together with symmetry conditions. Explicitly, the first four eigenfunctions are
\begin{eqnarray}
  && \tex{
   k=1: \quad \alpha_1=\frac{1}{2}, \quad f_1(y) = \frac{1}{2} y^2
   } ,  \nonumber \\
   &&
   \tex{
    k=2: \quad \alpha_2=1, \quad f_2(y) = 1 + \frac{1}{8N(N+2)} y^4
    },  \nonumber \\
   && \tex{  k=3: \quad \alpha_3=\frac{3}{2}, \quad f_3(y) = \frac{1}{2} y^2 + \frac{1}{48(N+2)(N+4)}y^6
   } ,  \nonumber \\
   && \tex{
    k=4: \quad \alpha_4=2, \quad f_4(y) = 1 + \frac{1}{4N(N+2)} y^4 + \frac{1}{192N(N+2)(N+4)(N+6)} y^8,
    }  \nonumber
\end{eqnarray}
 where the above properties are immediately apparent.

\begin{figure}[htp]
\vskip -2.2cm
\hspace*{-1cm}
\includegraphics[scale=0.6]{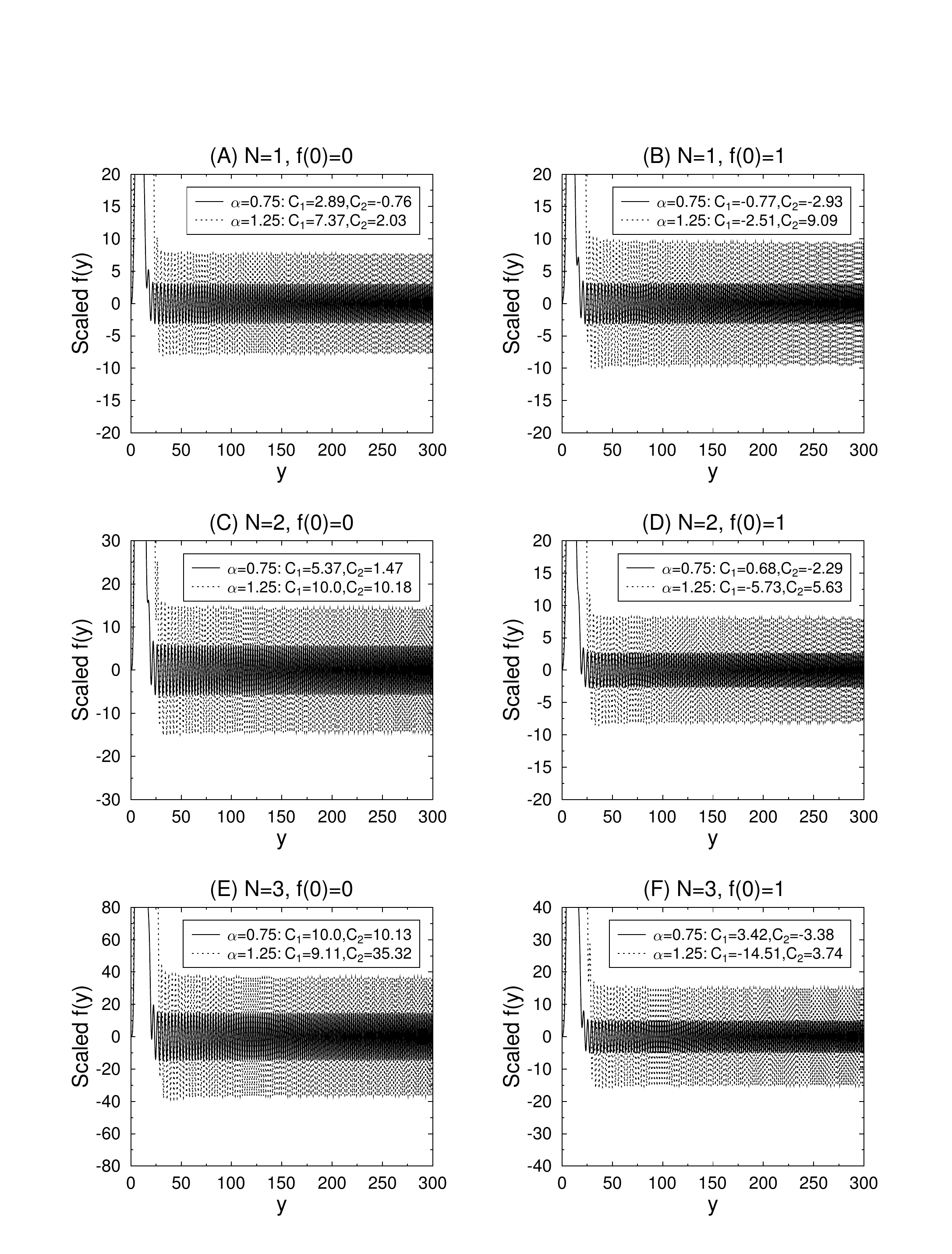}
\vskip -0.5cm \caption{ \small Numerical illustration of the oscillatory component in the far-field behaviour (\ref{b51}) in the linear case.
In each dimensional case $N=1,2,3$ scaled profiles (\ref{scaledf}) are shown for two selected values of the parameter $\alpha$
for each initial condition (\ref{sh1}) and (\ref{sh2}). The extracted least squares values of the far-field constants $C_{1,2}$ are stated inset in each figure.   }
 \label{fig1}
\end{figure}

\begin{figure}[htp]
\vskip -2.2cm
\hspace*{-1cm}
\includegraphics[scale=0.6]{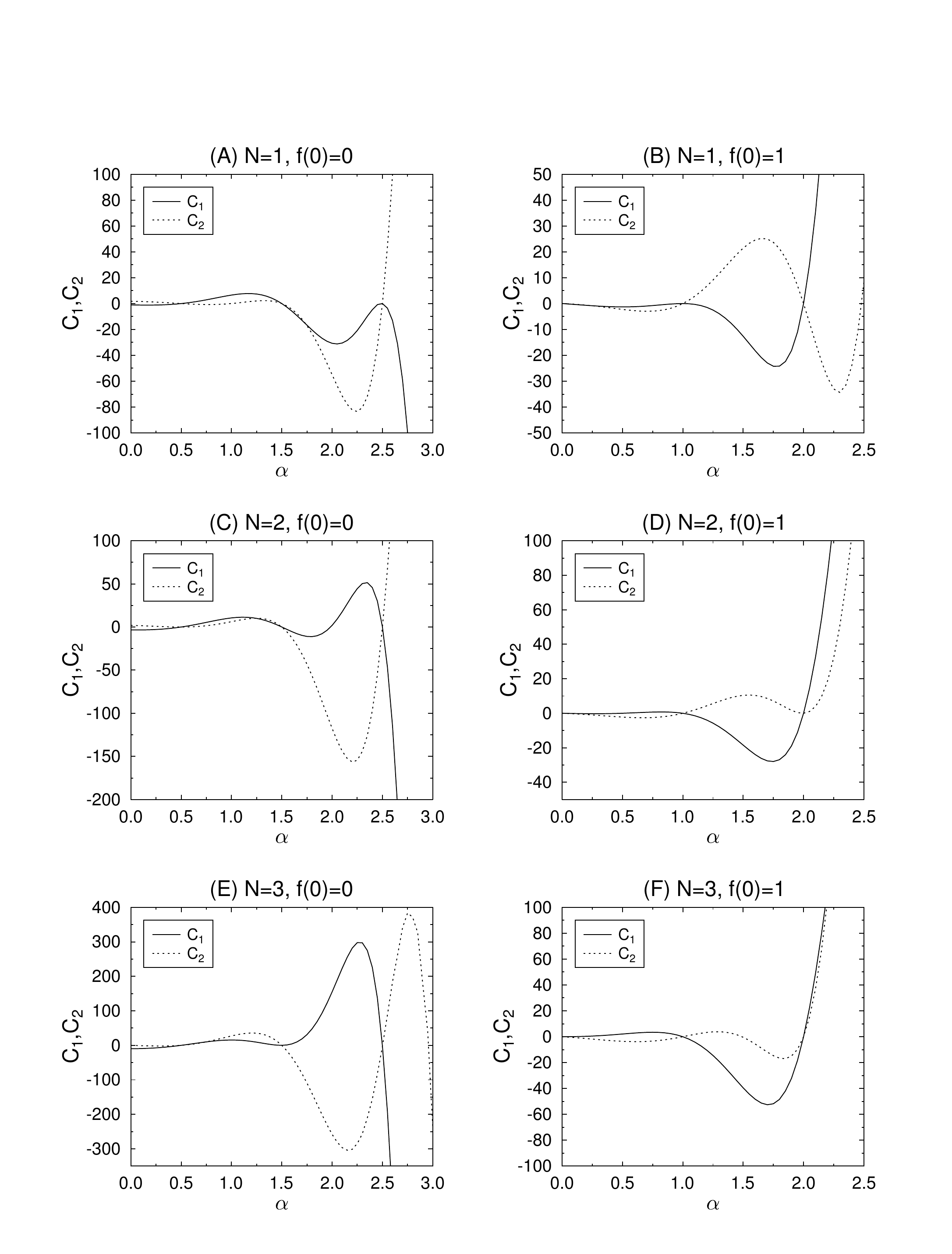}

\vskip -0.5cm \caption{ \small Numerical determination of the far-field constants $C_{1,2}$ for (\ref{b51}) in the stated
dimensional and initial condition cases. Shown are their variation with the parameter $\alpha$. The coincident zeros
correspond to the vanishing of the far-field maximal bundle yielding the eigenvalues of the linear problem $\alpha=\alpha_k$.  }
 \label{fig2}
\end{figure}

 Numerically we may illustrate the appearance of the eigenfunctions and their eigenvalues through consideration of an Initial Value Problem (IVP).
 The linear ODE (\ref{b2}) is solved numerically using ode15s (with tight error tolerances AbsTol=RelTol=$10^{-13}$) and subject to either (\ref{sh1})
 or (\ref{sh2}) as initial conditions. The far-field behaviour (\ref{b51}) is extracted from the numerical solution. The oscillatory component in the far-field behaviour
 is revealed by considering the scaled solution
 \be
\label{scaledf}
      f(y) y^{\frac{2}{3}(N+2\alpha)} e^{-c_0 y^{4/3}}.
 \ee
A least squares fitting of this function to the remaining oscillatory component in (\ref{b51}) over a suitable interval for large $y$, allows the
determination of the constants $C_{1,2}$. The large $y$ interval taken was typically [250,300]. The scaled function (\ref{scaledf}) is shown
in Figure \ref{fig1} for each of the two types of initial conditions in the $N=1,2$ and $3$ cases. Two selected values of the parameter $\alpha$
are taken in each case. The figures illustrate convergence to the oscillatory part of the far-field behaviour, with the extracted least squares estimates of the
constants $C_{1,2}$ shown inset for each of the two $\alpha$ values. Figure \ref{fig2} shows the variation of the far-field constants $C_{1,2}$ with
$\alpha$ in the $N=1,2,3$ cases and the two initial conditions. The constants $C_{1,2}$ both vanish precisely at the eigenvalues,
when $\alpha=\alpha_k$, the first five being shown in Figure \ref{fig2} in each $N$ case. The magnitude of the constants $C_{1,2}$
grow rapidly as $\alpha$ increases, making determination of further eigenvalues more difficult.

This approach demonstrates the validity of numerical determination of the eigenvalues and eigenfunctions by choosing $\alpha$ to minimise the far-field behaviour (\ref{b51}).

\section{A ``homotopic" transition to small $n>0$: some key
issues}
\label{S5}

    Next, we are going to use the above linear
    results to predict true nonlinear eigenvalues $\a_k(n)$ for
    the TFE--4, at least, for sufficiently small $n>0$. By
    continuity (to be discussed later on  as a continuous deformation after applying a homotopic argument) we now know that
     \be
     \label{b10}
      \tex{
     \a_k(0)= \frac k2, \quad k=1,2,3,...\,.
      }
      \ee
Another important conclusion: since, for $n=0$, any $f_k(0) \not
 = 0$ (this could happen only accidentally, with a probability 0),
 we can also expect that
  \be
  \label{b11}
  \mbox{for small $n>0$,} \quad f_k(0) \not = 0,
  \ee
  meaning that, in this case, we do not need to perform a shooting
  from the interface point, with a quite tricky behaviour nearby.

Also, we conclude from \ef{b51} that:
 \be
 \label{b12}
 \mbox{for small $n>0$, maximal bundle is 2D and oscillatory as $y
 \to + \iy$}.
  \ee
Those properties will aid progress on the nonlinear eigenvalue
problem.


 \subsection{Passing to the limit $n \to 0^+$ in the nonlinear
   eigenvalue problem.}


   Subsequently, we perform a homotopy deformation from the self-similar equation $\eqref{rad22}_{-}$ (with lower signs) of the TFE-4 \eqref{i1}
to the linear radial ODE \eqref{b2} corresponding to
the \emph{classical bi-harmonic parabolic equation} \eqref{bi.1}. In particular we construct
a continuous deformation from the radial equation $\eqref{rad22}_{-}$ to the linear equation \eqref{b2} for which we
know solutions explicitly.

It is clear that the CP for the \emph{bi-harmonic equation} \eqref{bi.1}
is well-posed and has a unique solution given by the convolution
 \be
 \label{b.11}
 u(x,t)=b(x-\cdot,t)\, * \, u_0(\cdot),
  \ee
   where $b(x,t)$ is the fundamental solution of the operator $D_t + \Delta^2$.

   Due to our analysis is possible to establish a connection between the radial solutions of \eqref{i1} and the corresponding to \eqref{bi.1}
   via the self-similar associated equations when $n \to 0^+$. To this end we apply the Lyapunov--Schmidt method to ascertain
   qualitative properties of the self-similar equation $\eqref{rad22}_{-}$ following a similar analysis as the one carried out in \cite{TFE4PV}.

   Thus, as we already know, the operator ${\bf A}^-_{0,y} (\a,f)$ defined by \eqref{b2} produces a countable family of eigenvalues
   \[
 \tex{
   \a_k\equiv \a_k(0)=\frac{k}{2}, \quad \hbox{with}\quad k=0,1,\cdots \, .
 }
   \]
      Note also that, \ef{b2} admits a complete and closed set of
 eigenfunctions being {\em generalised Hermite polynomials}, which
 exhibit finite oscillatory properties.

 This oscillatory issue seems to be crucial. In fact, in \cite{EGK1} was observe that a similar analysis of blow-up patterns
 for a TFE-4 like \eqref{i1} did not detect any stable
oscillatory behaviour of solutions  near the interfaces of the radially symmetric
associated equation. Hence, all the blow-up patterns turned out to be
nonnegative, which is a specific feature of the PDE under
consideration therein. However, this does not mean that blow-up similarity solutions of the CP do not change sign near the interfaces
 or inside the support. Actually, it was  pointed out
   that local sign-preserving property could be attributed only to the blow-up ODE and not to
   the whole PDE \eqref{i1}. Hence, the possibility of
having oscillatory solutions cannot be ruled out for every case.


 \subsection{Branching/bifurcation analysis}

  Now, we construct a continuous deformation such that the patterns occurring
 for the nonlinear eigenvalue problem $\eqref{rad22}_{-}$ are homotopically connected to the ones of the equation \eqref{b2}. Thus,
 we assume
for small $n>0$ in \eqref{rad22} the following expansions:
\begin{equation}
\label{br3}
    \a_k(n):= \a_k+ \mu_{1,k} n+ o(n),\quad
    |f|^n \equiv |f|^n= {\mathrm e}^{n\ln |f|}:=
     1 +n \ln |f|+o(n),
\end{equation}
 where the last one is assumed to be understood in a weak sense. The
second expansion cannot be interpreted pointwise for oscillatory
changing sign solutions $f(y)$, though now these functions are
assumed to have {\em finite} number of zero surfaces (as the
generalised Hermite polynomials for $n=0$ do). Indeed, as discussed in \cite{TFE4PV} for \eqref{rad22}
this is true if the zeros are transversal.

Furthermore, in order to apply the Lyapunov--Schmidt branching analysis we suppose the expansion
\begin{equation}
\label{br12}
 \tex{
    f=\sum_{|\b|=k} c_\b f_\b +V_k, \quad \hbox{for every $k\geq 1$,}
 }
\end{equation}
under the natural ``normalising" constraint
\begin{equation}
\label{nor}
 \tex{
    \sum\limits_{|\b|=k} c_\b=1.
    }
\end{equation}
Moreover, we write
 \[\{f_\b\}_{|\b|=k}=\{f_1,...,f_{M_k}\},\]
 as the natural basis of the $M_k$-dimensional eigenspace, with $M_k\geq 1$, such that
 \[
 \tex{
 f_k = \sum_{|\b|=k} c_\b f_\b \quad \hbox{and}\quad  V_k=\sum_{|\b|>k} c_\b f_\b,
 }
 \]
 with $V_k \in Y_k$ where $Y_k$ is the complementary invariant subspace of corresponding kernel. In particular, we
 consider the expansion
 \[ V_k:=n \Phi_{1,k} + o(n).\]
 Thus, applying the Fredholm alternative, we obtain the existence of a number of branches emanating
 from the solutions $(\a_k,f_k)$ at the value of the parameter $n=0$. We can guarantee that the first profile is unique but for the rest
 there could be more that branch of solutions emanating at $n=0$.
 The number of branches will depend on the dimension of the eigenspace so that the dimension is somehow involved; see \cite{TFE4PV} for any further comments and a detailed
 analysis of a similar branching analysis.

  \vspace{0.2cm}

\noindent{\bf Remark.} Furthermore, when $n\to 0^+$ we have a few different profiles of $f(y)$. For the PME--2 \eqref{pormed} the Graveleau profiles
are always unique by the Maximum Principle but not for the TFE--4 \eqref{i1}.
However, for $\a=\frac{1}{2}$ the only profile is $f(y)=y^2$.

   \vspace{0.2cm}

The previous discussion can be summarised in the following Lemma.
\begin{lemma}
The patterns occurring for equation \eqref{b2}, i.e.  the ones
for $\eqref{rad22}_{-}$ for $n=0$, act as branching points of
nonlinear eigenfunctions of the Cauchy problem for the operator
$\eqref{rad22}_{-}$, at least when the parameter $n$ is
sufficiently close to zero.
\end{lemma}

 \vspace{0.2cm}

\noindent{\bf Remark.} It turns out that using classical branching theory
 ``nonlinear eigenfunctions" $f(y)$ of changing sign, which
satisfies the {\em nonlinear eigenvalue problem} $\eqref{rad22}_{-}$ (with
an extra ``radiation-minimal-like" condition at infinity)
 at least, for
sufficiently small $n > 0$, can be connected with eigenfunctions $f_k$ of the linear problem \eqref{b2}.

 \section{Non-improvable regularity for the TFE--4}

The following main result is a straightforward consequence of our
focusing self-similar analysis.

\begin{theorem}
\label{theneg2}
 Let, for a fixed $n>0$, the nonlinear eigenvalue problem
 $(\ref{rad22})_-$
 have a nontrivial solution (eigenfunction) $f_k(y)$
 for some eigenvalue $\a_k>0$, i.e.  there exists a self-similar
 focusing solution of the problem $(\ref{eigpm})_-$ exhibiting the finite-time trace
 $\ef{gg5}$. Then:
\begin{enumerate}
 \item[(i)] For the general Cauchy problem for the TFE--4 $\ef{i1}$,
 even in the radial setting, the H\"older continuity  exponent of
 solutions
cannot exceed\footnote{We mean $C^{l+\e}$, if $\mu_k \ge 1$, that,
not that surprisingly,  happens for all small $n>0$.}
 \be
 \label{h11}
  \tex{
 \mu_k(n,N) = \frac {\a_k}{ \b_k} \equiv \frac {4 \a_k}{1+ n \a_k}.
  }
   \ee

  \item[(ii)] Let there exist $\mu_k<1$, i.e.
$\a_k < \frac 1{4-n}$, for $n \in (0,2)$\footnote{Actually, for
smaller $n$'s; Note that for every larger $n$'s solutions of the TFE--4 are known
to be strictly positive, \cite{BF1}).}. Then, for the TFE--4
\ef{i1},
 \be
 \label{h12}
  \tex{
  \n u_*(r,0^-) \in L^p_{\rm loc}(\ren) \quad \mbox{iff} \quad
  p \in[1,p_*), \,\,\, p_*(n,N)= \frac N{1- \mu_k},
  }
  \ee
 so that, for  (even radial) solutions of \ef{i1}, in general,
 for any $t>0$,
  \be
  \label{h13}
   \n u(x,t) \not \in L^{p_*}_{\rm loc}(\ren).
    \ee
    \end{enumerate}
 \end{theorem}

 \vspace{0.2cm}

 \noindent{\bf Remark.}  Basically $|\nabla u_*(r,0^-)| \in L^p_{\rm loc}(\ren)$ if and only if
   $$\tex{\int_{\re^N} |\nabla u_*(r,0^-)|^p = C \int_0^1 r^{N-1} (r^{\mu_k-1})^p <\infty.}$$
   Thus, since $0<\mu_*<1$ we find that to have that integral bounded we need that, evaluating the coefficients in the second integral
   $$\tex{ p <p_*\equiv p_*(n,N) =\frac{N}{1-\mu_k(n,N)}.}$$

  \vspace{0.2cm}

\noindent{\bf Remark.} Of course, it can happen that $\mu_k>1$, so focusing does not
supply us with a truly H\"older continuous focusing trace but
rather than a $C^{l+\e}$-trace. Actually, by continuity in $n$,
exactly this happens for $n>0$ small, where $\mu_1(0,N)=2$.

We expect that the minimal value of $\mu_k$ in \ef{h11} is
attained at $k=1$, but cannot prove this for all $n>0$ (for small
ones, this is obvious). Note again that, for all large $n \ge 2$,
such a focusing is not possible in principle, since the solutions
must remain strictly positive for all times, \cite{BF1}.

Thus, this optimal (non-improvable) regularity results for the
TFE--4 (and, seems for many other parabolic equations) depends on
the solvability of the nonlinear eigenfunction focusing problem.

Note that in the
previous section we showed a homotopy connection between that nonlinear eigenvalue
problem and the linear problem at $n=0$ \eqref{b2}. For those reasons
the analysis performed in those two previous sections is crucial
in ascertaining qualitative information
about the solutions of the TFE--4.

Furthermore, we can also ascertain the H\"{o}lder continuity with respect to the temporal variable $t$. This fact comes directly from 
the non-improvable H\"{o}lder's exponent we have already obtained above. Indeed, assuming the radial self-similar solutions of the form \eqref{sssol} 
for the focusing, 
one easily finds that 
$$|u(0,t)-u(0,0)|=(-t)^\a f(0).$$  
Hence, the $t$-H\"{o}lder's exponent close to $t=0$ cannot be bigger than $\a$. 

Therefore, the focusing solution given here provides us with 
optimal H\"{o}lder estimations for the variables $x$ and $t$.

 \section{Oscillatory structure of maximal solutions for $n>0$}
  \label{S.maxosc}

\begin{figure}[htp]
\vskip -.2cm \hspace*{-1cm}
\includegraphics[scale=0.6]{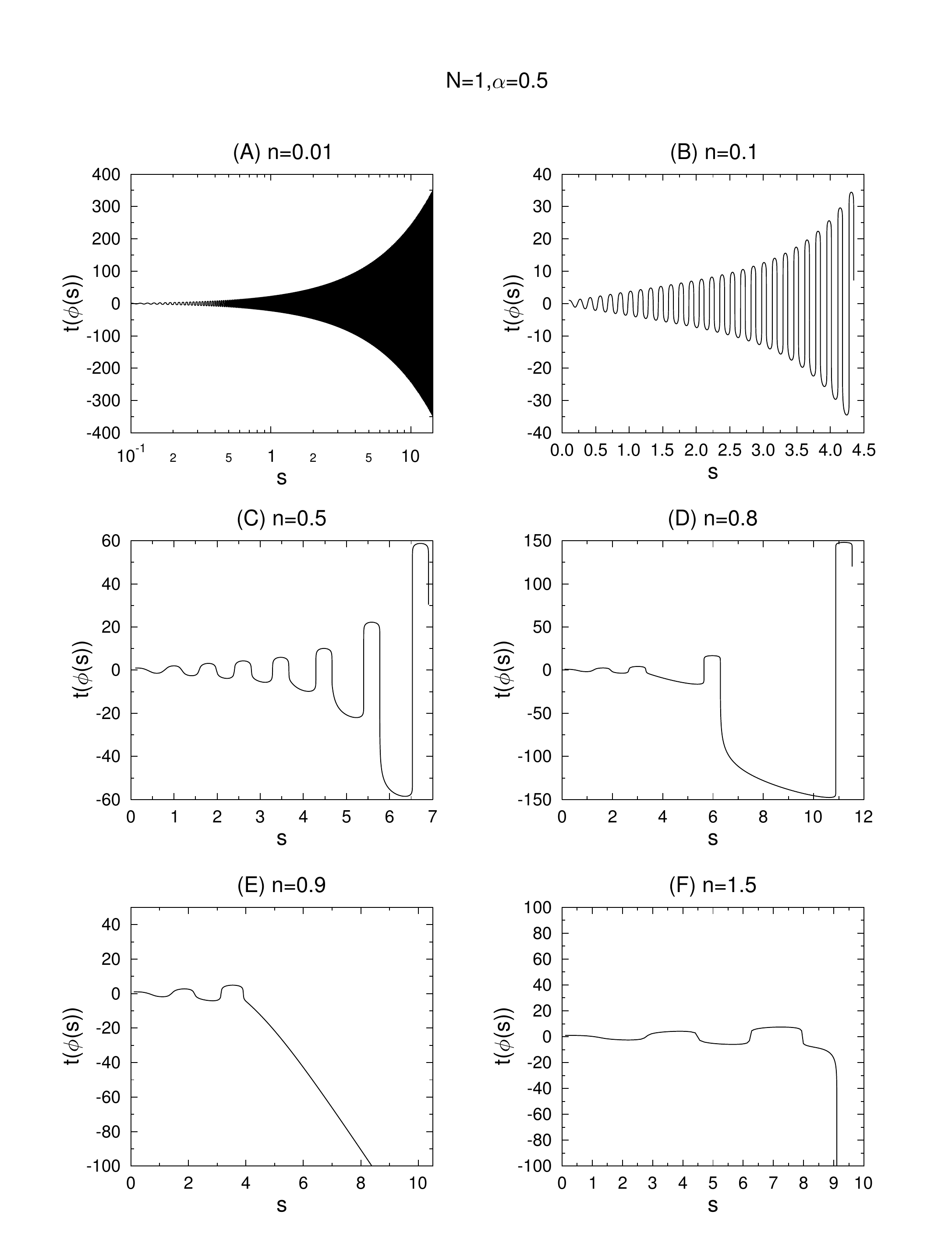}

\vskip -0.5cm \caption{ \small Illustrative numerical profiles of transformed $\phi(s)$ via (\ref{tphi}).
Shown is the parameter case $N=1,\alpha=0.5$ for selected $n$.
   }
 \label{fig3}
\end{figure}

Let us now consider $n \in (0,2)$. For the ODE \ef{rad22}, we will try the
same {\em anzatz} as in \cite{EGK2}. Namely, we will find
solutions of the maximal type, with the envelope as in  \ef{gg2}.
We introduce a corresponding {\em oscillatory} component as
follows:
 \be
 \label{zz1}
 f(y)= y^{\mu} \var(s), \quad s= \ln y, \quad \mu = \frac{4}{n} \quad \mbox{for} \quad
 y \gg 1.
 \ee
 Then,
  \[
  \tex{
  f'=y^{\mu -1} \big( \dot{\var} + \mu \var \big), \quad
  f''=y^{\mu -2} \big[ \ddot{\var} + \big( \mu - 1 \big) \dot{\var} + \mu \big(\mu-1 \big) \var \big],
  }
  \]
  \[
    f''' = y^{\mu-3} \left[ \dddot{\var} + 3(\mu-1) \ddot{\var} + (3 \mu^2 - 6\mu +2) \dot{\var}
     + \mu(\mu^2 -3 \mu + 2)\ \var
       \right],
  \]
  where $'={\mathrm d}/{\mathrm d}y$ and $\dot{}={\mathrm d}/{\mathrm d}s$. Substituting these expressions into \ef{rad22} yields the following
  fourth-order homogeneous ODE for $\var(s)$
\begin{align}
 &  \ddddot{\var} + 2(N-4 + 2\mu) \dddot{\var} + (6 \mu^2 + 6(N-4)\mu + 11 + (N-1)(N-9)) \ddot{\var}  \label{phieqn}\\
 &    + 2(2\mu+N-4)(\mu^2 + (N-4)\mu +2-N) \dot{\var}
  + \mu (\mu-2) (\mu^2 +2(N-3)\mu + 3+ (N-1)(N-5)) \var  \nonumber \\
 &   \tex{
 + n \left( \frac{\dot{\var}}{\var} + \mu  \right) \left[ \dddot{\var} + (N-4+3\mu) \ddot{\var}
     + (3\mu^2 + 2(N-4)\mu +4-2N) \dot{\var} + \mu (\mu-2)(N-2+\mu) \var \right]
     }
       \nonumber \\
 &   \tex{
 + \left(  \frac{1}{4}(1+n\alpha)(\dot{\var} + \mu \var) -\alpha \var \right) |\var|^{-n} = 0 .
 } \nonumber
\end{align}
We mention that this equation is autonomous and thus may be reduced to third-order, although we do not utilise
this reduction here. Figure \ref{fig3} illustrates the periodic nature of the solutions for $\phi$, at least for $n$ small enough.
Since the oscillations occur over such a large range, we use the following transformation to allow the oscillations to be visible on the plots,
\be
    \mbox{t}(\phi(s)) = \left\{ \begin{array}{ll}
                          \ln \phi(s) +1, & \mbox{if}\; \phi(s)>1, \\
                          \phi(s), & \mbox{if}\; -1<\phi(s)<1, \\
                         - \ln (-\phi(s)) -1, & \mbox{if}\; \phi(s)<-1.
                          \end{array}\right.
\label{tphi}
\ee
Plotted are numerical solutions for $\phi$ in the two cases $\alpha=0.5,1$ with $N=1$. The subplots illustrate the change
in the profile behaviour with $n$. The numerics support the proposition that the family
\ef{zz1} is 2D, composed of:
\begin{enumerate}
\item[(i)] a 1D stable manifold $\var_*(s)$,\\
and
\item[(ii)] a phase shift $\var_*(s+s_0)$ for any $s_0 \in \re$.
\end{enumerate}

As we have seen earlier, these two properties are true for $n=0$,
and, hence, by continuity, we conjecture remain true for small $n>0$.

However,
the periodic exponential structure, with the linear ODE for $n=0$ is
replaced by a more difficult one \ef{zz1} for $n>0$. Nevertheless, the
periodic nature of such a behaviour appears universal, and
is expected to remain for larger $n$.
The numerics though do suggest that it is lost
(in a homoclinic-heteroclinic bifurcation) for $n$ near to 0.8 or 0.9.

\begin{figure}[htp]
\vskip -2.2cm
\hspace*{-1cm}
\includegraphics[scale=0.6]{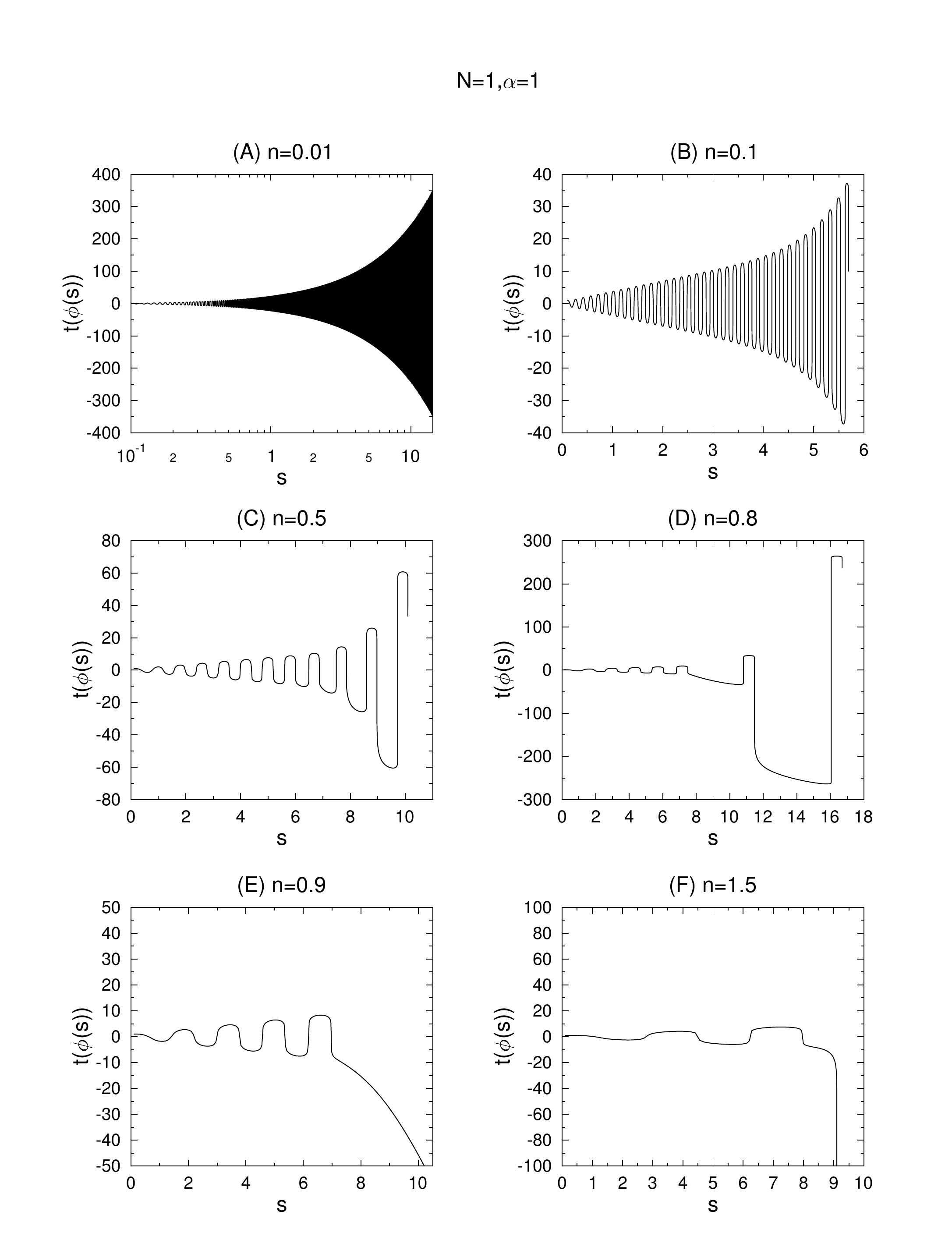}

\vskip -0.5cm \caption{ \small Illustrative numerical profiles of transformed $\phi(s)$ via (\ref{tphi}).
Shown is the parameter case $N=1,\alpha=1$ for selected $n$.
   }
 \label{fig4}
\end{figure}

To reconcile with the known behaviour in the linear case,
we now study the behaviour of the periodic solutions for small $n>0$. To reveal the limiting oscillatory
behaviour as $n \to 0$, we keep the leading terms of the coefficients in (\ref{phieqn}) as $\mu \to \infty$ to obtain
\be
   \ddddot{\var} + 4\mu \dddot{\var} + 6 \mu^2  \ddot{\var} + 4\mu^3 \dot{\var} + \mu^4 \var
 + \frac{1}{4} \left(  \dot{\var} + \mu \var \right) |\var|^{-n} = 0 . \label{phismalln}
\ee
We now rescale as follows
 \be
 \label{AAA.1}
 \tex{
\var(s)=  \mu^{-\frac{3}{n}} \hat{\var}(\hat{s}), \quad \displaystyle s = \frac{\hat{s}}{\mu}   ,
  }
\ee
leading to
\be
 \tex{
   \ddddot{\hat{\var}} + 4 \dddot{\hat{\var}} + 6  \ddot{\hat{\var}} + 4 \dot{\hat{\var}} + \hat{\var}
 + \frac{1}{4} \left(  \dot{\hat{\var}} + \hat{\var} \right) |\hat{\var}|^{-n} = 0 ,
 }
  \label{phihat}
\ee where here $\dot{}$ denotes ${\mathrm d}/{\mathrm d}\hat{s}$.
For $n=0$, the equation (\ref{phihat}) becomes linear,
\be
  \tex{
    \ddddot{\hat{\var}} + 4 \dddot{\hat{\var}} + 6  \ddot{\hat{\var}} + 4 \dot{\hat{\var}} + \hat{\var}
 + \frac{1}{4} \left(  \dot{\hat{\var}} + \hat{\var} \right) = 0 ,
 }
  \label{phihatn0}
\ee
with the characteristic equation
$$
 \tex{
(m+1)\left( (m+1)^3 + \frac{1}{4} \right)=0,
 }
 $$
for exponential solutions
$
\hat{\var}(\hat{s})=  {\mathrm e}^{m \hat{s} } .
$
The roots of the characteristic equation are linked to those of the controlling factor in (\ref{b3}) being
$$
 \tex{
   m +1 = \frac{4}{ 3}a_1,\quad \frac{4}{3}a_2,\quad  \frac{4}{3}a_3 \quad \mbox{and} \quad 0 .
 }
 $$
Consequently we obtain the dominant asymptotic behaviour
\be
 \tex{  \hat{\var}(\hat{s}) \sim e^{(-1+\frac{4}{3}c_0) \hat{s} } \left[ \hat{A}_1 \cos \left(\frac{4}{3}c_1\hat{s} \right)
      + \hat{A}_2  \cos \left(\frac{4}{3}c_1\hat{s} \right) \right] \quad
       \mbox{as $\hat{s} \to \infty$},
 }
 \label{phihatasy}
\ee
for arbitrary constants $\hat{A}_{1,2}$. Denoting the independent variable by $\hat{y}$ rather than $y$ for convenience, we thus have the small $n$ behaviour
\be
 \tex{
   f(\hat{y}) \sim \hat{y}^{\frac{4}{n}} \left( \frac{4}{n}\right)^{-\frac{3}{n}} \hat{\var} \left({\frac{4}{n}\ln(\hat{y})} \right) .
 }
 \ee
This may be reconciled with the expression in (\ref{b51}) through the identifications
\[
  \tex{
      \hat{y} = e^{\frac{3n}{16} {y}^{4/3}}, \quad
     \hat{A}_{1,2} = \left( \frac{4}{n}\right)^{\frac{3}{n}} C_{1,2} .
 }
 \]

 \section{Nonlinear eigenfunctions by shooting}

 The eigenfunctions for $n>0$ may be obtained via shooting. The conditions (\ref{sh1}) or (\ref{sh2}) are again
 used as initial conditions and $\alpha$ determined by capturing the growth (\ref{mingro}) for sufficiently large $y$ values. This growth behaviour may be imposed by
 $$\hbox{minimising}\quad \mu y f' -\alpha f\quad \hbox{at}\quad y\quad \hbox{values},$$
 typically chosen to be around 40. Standard regularisation of the term $|f(y)|^n$ is required in the form $(f^2+\delta^2)^{\frac{n}{2}}$ with $\delta$ taken relatively small.

\begin{figure}[htp]
\hspace*{-1cm}
\includegraphics[scale=0.6]{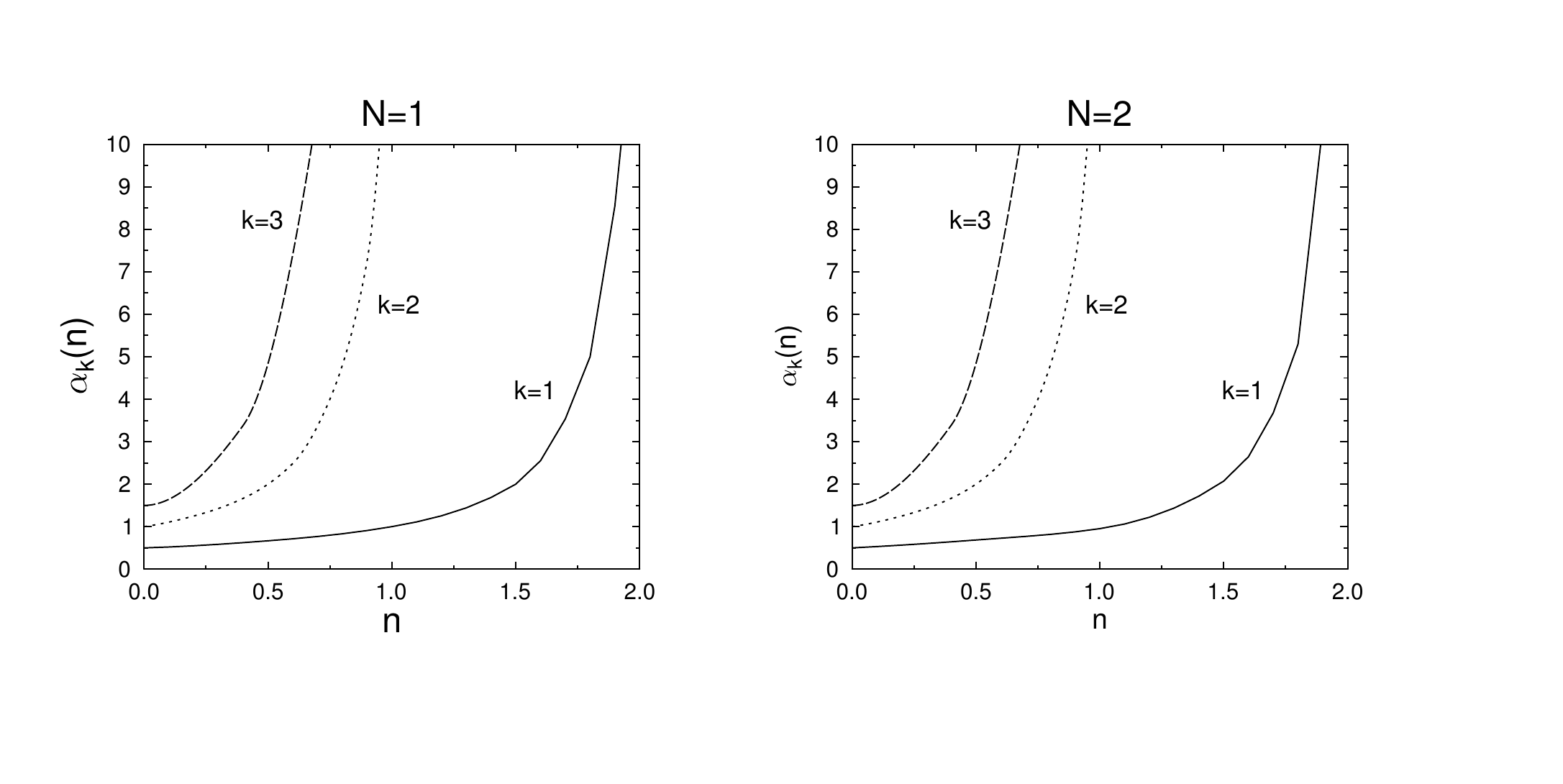}

\vskip -0.5cm \caption{ \small Numerical determination of the eigenvalues $\alpha_k(n)$ for the first three branches $k=1,2,3$ of eigenfunctions.
   }
 \label{fig5}
\end{figure}

 Figure \ref{fig5} shows the eigenvalues of the first three branches $k=1,2,3$. These numerically extend the $n=0$ eigenvalues of Proposition \ref{Pr.anal} to $n>0$.

 As in the $n=0$ linear case, the eigenvalues remain the same irrespective of the spatial dimension $N$, at least to the accuracy of the numerical calculations that were performed.

 Furthermore, it is worth remarking that along each branch i.e.  fixed $k$, the exponent $\mu$ of the far-field behaviour of the eigenfunctions
 $$f\sim C y^{\mu},$$
 remains quite close to its value in the linear case i.e.,
 $$\mu(n)\approx \mu(0).$$
 Thus, as an approximation we have
 $$\alpha_k(n) \approx \alpha_k(0)/(1-\alpha_k(0) n),$$
 which seems to be a reasonable approximation at least
 $$\hbox{for}\quad n< 1/\alpha_k(0).$$
 The corresponding eigenfunctions are thus similar to those in the linear case.

An alternative approach to determining the eigenfunctions in both the linear $n=0$ and nonlinear $n>0$ cases,
is through minimisation of the oscillatory maximal profile to obtain the non-oscillatory minimal profile for $y\gg 1$.

In principle we have two parameters $\alpha$ and, say $\nu=f''(0)$ or $f(0)$, depending on which branch of eigenfunctions is considered.
In the linear case, using \ef{b51}, we are required to satisfy two algebraic equations with analytic functions:
     \be
     \label{sh3}
     \left\{
     \begin{matrix}
     C_1(\a,\nu)=0, \ssk \\
     C_2(\a,\nu)=0.
      \end{matrix}
      \right.
      \ee
Therefore, we arrive at a well-posed $2-2$ shooting problem, which
cannot have more than a countable pairs of solutions (as mentioned already). This approach is practical for the linear
$n=0$ case, as the two parameter form of the maximal bundle is known explicitly as given in (\ref{b51}) and was essentially
pursued in section 4.2. However, lacking such an explicit expression for $n>0$, limits this approach in the nonlinear case.

\begin{figure}[htp]
\vskip -2.2cm
\hspace*{-1cm}
\includegraphics[scale=0.6]{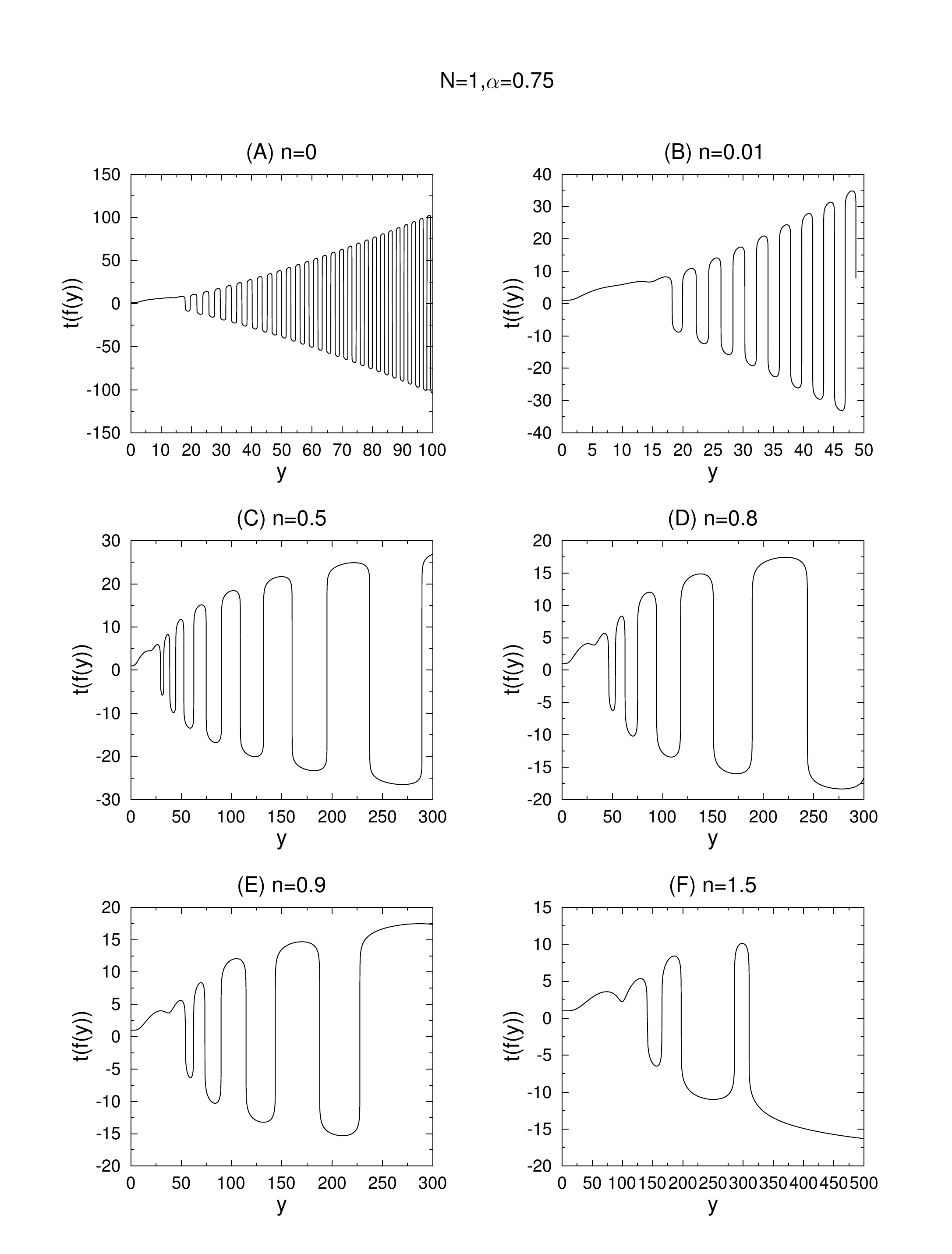}

\vskip -0.5cm \caption{ \small Illustrative numerical profiles of transformed $\phi(s)$ via (\ref{tphi}). Shown is the parameter case $N=1,\alpha=1$ for selected $n$.
   }
 \label{fig6}
\end{figure}

Finally it is worth presenting some numerical experiments showing the behaviour of the maximal profiles. Figure \ref{fig6}
shows profiles in the case $N=1,\alpha=0.75$ for selected $n$. Since the oscillatory profiles occur over such large ranges, the transformed profile using (\ref{tphi}) is depicted.

The figures suggest the loss of a pure oscillatory structure as $n$ increases, which is compatible with the behaviour
seen for $\phi$ in Figures \ref{fig3} and \ref{fig4}.

We conjecture that this is highly suggestive of a homoclinic-heteroclinic bifurcation. The precise values of $n$ though
at which this occurs is difficult to determine with sufficient accuracy.

\vspace{0.1cm}

 \section{After-focusing self-similar extension}

As usual and as we have observed in many similar problems on a
(unique) extension of a solution after blow-up, this is much
easier. In fact, after focusing at $t=0^-$, we arrive at a
well-posed CP (or TFE, that does not matter since $u_0(r)>0$, no
interfaces are available) for the TFE--4 \ef{i1} with initial data
\ef{gg5} (with $\mu=\mu_k$ when possible), already satisfying the
necessary minimal growth at infinity.

Therefore, there exists a self-similar solution $(\ref{upm})_+$
with $f(y)$ satisfying the ODE $(\ref{rad22})_+$, {\em with
already fixed ``eigenvalue" $\a=\a_k$}. Because of changing the
signs in front of the linear terms, this changes the dimension of the
stable/unstable manifolds as $y \to +\iy$ (in particular, the
unstable manifold becomes 1D, similar to $n=0$).

To explain the latter, as usual consider the case $n=0$ (and hence
small $n>0$). Then, with the change of sign in the two linear
terms in \ef{b2}, we arrive at a different characteristic equation
for $a$:
 \be
 \label{b3N}
  \tex{
  f(y) \sim \eee^{a y^\g} \asA y \to + \iy \LongA a^3= {\bf +} \frac 14 \,\big(
  \frac 34 \big)^3.
   }
   \ee
 This characteristic equation gives two stable roots with ${\rm Re}(\cdot)<0$:
  \be
  \label{b4N}
   \tex{
  b_{1,2}=\frac 34\,4^{-\frac 13} \big[ -\frac 12 \pm {\rm i} \, \frac{\sqrt 3}2
  \big] \equiv c_0 \pm {\rm i} \, c_1.
 }
 \ee
 and one positive root maximal solutions;
 see below)
  \[
   \tex{
 b_3= \frac 34 \, 4^{-\frac 13} > 0.
 }
 \]
 Hence, the bundle of maximal solutions is oscillatory as $y \to +
 \iy$  is just 1D, so we arrive at an {\em undetermined} problem
 for $\nu,\a$, which satisfy just a single algebraic equation (unlike two in \ef{sh3}).
 Indeed, this makes the $\a$-spectrum continuous and solvability
 for any $\a>0$.

\ssk

This allows us to get a unique extension profile $f(y)$ for such a
{\em a priori} fixed eigenvalue $\a_k$. In other words, the
focusing extension problem is not an eigenvalue one, since a
proper $\a_k$ has been fixed by the previous focusing blow-up
evolution.

It might be said, using standard terminology from linear operator
theory, that for the sign $+$ the spectrum of this problem
becomes continuous, unlike the discrete one for $-$. Therefore, we
do not study this much easier problem anymore, especially, since
it has nothing to do with our goal: to detect an non-improvable
regularity for the TFE--4 via blow-up focusing.


\appendix
\setcounter{equation}{0}
\section{The limit $n \to 0$ for maximal solutions}
\label{A1}

We consider here the singular limit $n\to 0^+$ for the
equation in (\ref{kk2}), written explicitly here as
\[
 \tex{
  \frac{1}{y^{N-1}} \left[ y^{N-1}|f|^n \left(\frac{1}{y^{N-1}}(y^{N-1} f')'\right)' \right]'
  + \frac{1}{4}(1+\alpha n)yf' - \alpha f=0 .
 }
 \]
The non-uniform
solution in this limit comprises two regions, an {\em Inner
region} where $1\ll y < O(n^{-3/4})$ with $|f|^n \sim 1$, and an {\em Outer region}
$y=O(n^{-3/4})$ where $\ln |f| = O(1/n)$.
The labelling of these regions as inner and outer becomes apparent
during the course of the scalings.

In the inner region $1\ll y < O(n^{-3/4})$, we obtain
at leading order in $n$ the linear ODE (\ref{b2}). Posing $f\sim f_0$ with $|f_0|^n \sim 1$, we obtain
\be
\label{inner1} \tex{
 \frac{1}{y^{N-1}} \left[ y^{N-1} \left(\frac{1}{y^{N-1}}(y^{N-1} f_0')'\right)' \right]'
  + \frac{1}{4}(1+\alpha n)yf_0' - \alpha f_0 =0.
 }
 \ee
The far-field behaviour of (\ref{inner1}) may be determined using a
WKBJ expansion in the form
\begin{equation}
   f_0(y) \sim A(y) {\mathrm e}^{\phi(y)} \hspace{1cm} \mbox{as $y \to +\infty$}
\label{inner2}
\end{equation}
which gives
\begin{equation}
     (\phi')^3 = \frac{1}{4}y, \hspace{1cm}
     3 y A'+ 2 \left( N+2\alpha-1 -12 (\phi')^2 \phi''\right) A =0.
\label{inner3}
\end{equation}
The required solutions to (\ref{inner3}) take the form
\begin{equation}
 \tex{
       \phi(y) = a y^{4/3} , \hspace{1cm}
             A(y) = k y^{-\frac{2}{3}(N + 2\alpha)} ,
             }
\label{inner4}
\end{equation}
where $a$ satisfies the cubic equation in (\ref{b3}) and $k$ an arbitrary constant. Thus, the dominant behaviour for large $y$ is
\be
 \label{inner5}
   f_0(y) \sim  y^{-\frac{2}{3}(N+2\alpha)} \left( k_1 \exp \big\{ a_1 y^{\frac{4}{3}} \big\}
   + k_2 \exp \big\{ a_2 y^{\frac{4}{3}} \big\} \right) ,
\ee with $a_{1,2}$ as given in (\ref{b4}) and $k_{1,2}$ complex
constants chosen so that  the expression is real as stated in
(\ref{b51}).

This inner solution breaks down when $y=O(n^{-3/4})$, where $\ln
|f_0| =O(1/n)$. This suggests the consideration of an outer region
with scaling $Y=n^{-\frac 3 4}y$. In $Y=O(1)$, we have
\be
 \label{outer1}
  \tex{
 \frac{n^3}{Y^{N-1}} \frac{\mathrm d}{{\mathrm d}Y}\left[ Y^{N-1}|f|^n
   \frac{\mathrm d}{{\mathrm d}Y} \left(\frac{1}{Y^{N-1}} \frac{\mathrm d}{{\mathrm d}Y}
    \left(Y^{N-1} \frac{{\mathrm d}f}{{\mathrm d}Y} \right)\right) \right]
  + \frac{1}{4}(1+\alpha n)Y \frac{{\mathrm d}f}{{\mathrm d}Y} - \alpha f=0 .
 }
 \ee
Rather than posing a multiple-scales ansatz directly, it is more convenient to work in complex form and
we may instead consider
\begin{equation}
      f(Y) \sim {\mathrm e}^{b(Y)/n} B(Y) \hspace{1cm} \mbox{as \,\,$n \to 0$},
\label{outer2}
\end{equation}
where $b$ is complex in order to match with the inner solution expression (\ref{inner5}).
Thus, at $O(1/n)$ in (\ref{outer1}) we obtain
\begin{equation}
 \tex{
 |{\mathrm e}^b| \left( b' \right)^3 + \frac{1}{4}Y = 0,
 }
 \label{outer3}
\end{equation}
whilst at $O(1)$ we have
\begin{equation}
 \tex{
 3\frac{B'}{B} +  \left( 1 - \frac{\alpha}{4} + \ln |B|\right)b' + 6\frac{b''}{b'} + \frac{2((N-1)+2\alpha)}{Y}
  = 0 ,
  }
 \label{outer4}
\end{equation}
where $'$ denotes $\frac {\mathrm d}{{\mathrm d}Y}$ and the approximation
\begin{equation}
 |f|^n \sim |{\mathrm e}^b| \left( 1 + n \ln |B|\right)
\label{outer5}
\end{equation}
has been used. The solution to (\ref{outer3}) that matches with (\ref{inner2}) and (\ref{inner4}) is
\begin{equation}
  \tex{
   b(Y) = a \frac{3}{c_0} \ln \left( 1 + \frac{c_0}{3} Y^{4/3} \right).
 }
 \label{outer6}
\end{equation}
The amplitude $B(Y)$ is determined from (\ref{outer4}) using the solution (\ref{outer6}). We obtain
\[
 \tex{
   \ln |B| = \frac{1}{\left(1+ \frac{c_0}{3} Y^{4/3}\right)} \left( |k| - \frac{2}{3}(N+2\alpha) \ln Y - \frac{c_0}{12}(2N+3\alpha-4) Y^{4/3}  \right),
 }
 \]
after matching to (\ref{inner2}) with (\ref{inner4}).

\end{document}